\documentclass{article}%
\usepackage{amsfonts,float}
\usepackage{mathrsfs}
\usepackage{graphicx}
\usepackage{amsmath}
\usepackage{amssymb}%
\providecommand{\U}[1]{\protect\rule{.1in}{.1in}}
\newtheorem {theorem}{Theorem}[section]
\newtheorem {proposition}{Proposition}[section]

\newtheorem{lemma}{Lemma}[section]

\newtheorem{remark}{Remark}[section]

\newtheorem{assumption}{Assumption}
\newenvironment{proof}[1][Proof]{\textbf{#1.} }{\
\rule{0.5em}{0.5em}}

\usepackage{subfigure}

\begin{document}

\begin{center}
{\LARGE Variable bandwidth kernel regression estimation}



\bigskip Janet Nakarmi$^{a}$, Hailin Sang$^{b}$ $^{1}$\footnotetext[1]{Corresponding author} and Lin Ge$^{c}$

\bigskip$^{a}$ Department of Mathematics, University of Central Arkansas, Conway, AR 72035, USA. E-mail address: janetn@uca.edu

\bigskip$^{b}$ Department of Mathematics, The University of Mississippi,
University, MS 38677, USA. E-mail address: sang@olemiss.edu

\bigskip$^{c}$ Division of Arts and Sciences, Mississippi State University at Meridian,
Meridian, MS 39307, USA. E-mail address: lge@meridian.msstate.edu
\end{center}

\bigskip\textbf{Abstract}
In this paper we propose a variable bandwidth kernel regression estimator for $i.i.d.$ observations in $\mathbb{R}^2$ to improve the classical Nadaraya-Watson estimator. The bias is improved to the order of $O(h_n^4)$ under the condition that the fifth order derivative of the density function and the sixth order derivative of the regression function are bounded and continuous. We also establish the central limit theorems for the proposed ideal and true variable kernel regression estimators. The simulation
study confirms our results and demonstrates the advantage of the variable bandwidth kernel method over the classical kernel method.




\section {Introduction}


Let $(X_1, Y_1), ..., (X_n, Y_n)$ be $i.i.d.$ observations in $\mathbb{R}^2$ such that 
\begin{eqnarray*}\label {regfun}
Y_i=r(X_i) +\varepsilon_i
\end{eqnarray*}
where the $\varepsilon_1,..., \varepsilon_n$ are $i.i.d.$ random variables with $E \varepsilon_1=0$ and $E \varepsilon_1^2<\infty$, and $\varepsilon_i$ and $X_i$ are independent for every $i \in [1,n]$. $f(t)$ is the probability density function of $X_1$. One classical nonparametric estimator for the regression function $r(t)$ introduced independently by Nadaraya \cite{EAN1964}  and Watson \cite{GSW1964} is

\begin{eqnarray}\label {hatrt}
\hat{r}(t; h_n)=\frac{\hat{g}(t; h_n)}{\hat{f}(t; h_n)},
\end{eqnarray}

where 

\begin{eqnarray} 
\hat{g}(t; h_n)= \frac{1}{n h_n}\sum_{i=1}^n K\left ( \frac{t-X_i}{h_n}    \right )     Y_i , \nonumber\\
\hat{f}(t;h_n)= \frac{1}{n h_n}\sum_{i=1}^n K\left ( \frac{t-X_i}{h_n}    \right ). \label {hatft}
\end{eqnarray}

The kernel $K$ satisfies $K(x)\geq 0$ and $\int_\mathbb{R} K(x)dx=1$, and the bandwidth sequence $h_n$ satisfies $h_n\rightarrow 0$ and $nh_n\rightarrow \infty$ as $n \rightarrow \infty$. (\ref{hatft}) is the Parzen-Rosenblatt estimator for the density function $f(t)$ of random variable $X_1$. 
If $K$ is symmetric to zero, $f(t)$  and $r(t)$ have bounded second order derivatives, then the bias of (\ref{hatrt}) has order of $O(h_n^{2})$ and the variance has order of $O((nh_n)^{-1})$.  See, e.g., \cite{R}, \cite{GC1977}.
Noda \cite{KN1976} established the convergence of the Nadaraya-Watson estimator to $r(t)$ and the mean square error at a fixed point where $r(t)$ is continuous. 
The uniform consistency of the Nadaraya-Watson estimator was shown in \cite{HJB1983} for the case of discrete $X_i's$. See \cite{wasserman} and \cite{MS2000} and references therein for more literature on 
 the Nadaraya-Watson estimator. 
 
In general, for a measurable function $l$, the regression function $r(t,l):=\mathbb{E}(l(Y)|X=t)$ if it exists. Notice that for $l(y)=y$, we have the regression function $r(t)$. The Nadaraya-Watson type kernel regression estimator has form  
    \begin{equation}\label{rtq}
    \hat{r}_n(t,l)=\frac{\sum_{i=1}^n K\left(\frac{t-X_i}{h_n}\right)l(Y_i)}{nh_n\hat{f}(t;h_n)}.
    \end{equation} 
This estimator was studied by many authors, for example, Einmahl and Mason (\cite{EM2000}, \cite{EM2005}) obtained the exact rate of uniform consistency $O_{a.s.}(\sqrt{nh_n|\log h_n|})$ of regression estimator (\ref{rtq}) with some additional smoothness conditions for $l(\cdot)$ in a compact interval. 

For fixed value $t$ with $f(t)\ne 0$, in the Nadaraya-Watson regression estimation (\ref{hatrt}) for $r(t)$ or the Parzen-Rosenblatt density estimation for $f(t)$, we use the same bandwidth $h_n$ without considering the location of each data $X_i$ relative to $t$. The application of constant bandwidth everywhere makes bias to have the order of $h_n^2$.  To reduce the order of the bias, in kernel density estimation,  
Abramson \cite{Abramson}  applied the so called `square root law' which allows the bandwidth to vary with the data.  That is, if one takes the bandwidth $h_n/\gamma^{1/2}_t(X_i)$ in the classical density estimator,  with the same sequence $h_n\to 0$, where $ \gamma_t(x)=f(x)\vee (f(t)/10)$, then the estimator turns to be 
\begin{equation}\label{abr1}
f_n(t)=\frac{1}{n}\sum_{i=1}^n\frac{ 1}{h_n/\gamma^{1/2}_t(X_i)}K\left(\frac{t-X_i}{h_n/\gamma^{1/2}_t(X_i)}\right).
\end{equation}
The  bias of this estimator is reduced to the order of $h_n^4$ under the assumption that $f(t)\ne 0$ and $f$ has fourth order bounded and continuous derivatives around $t$. The true bandwidth $ h_n/\gamma_t ^{1/2}(X_i)$, at each observation $X_i$ is inversely proportional to $f^{1/2}(X_i)$ if $f(X_i)\ge f(t)/10$ (which is the square root law).  

It was once believed that $\gamma_t(x)$ in (\ref{abr1}) could be replaced by $f(x)$ and only the square root law with the bandwidth $h_n/f^{1/2}(X_i)$ makes the bias reduction work. But Terrell and Scott \cite{TerrellScott}  showed that the bias reduction can not be reached in some cases if only the square root law with this bandwidth $h_n/f^{1/2}(X_i)$ is applied. In fact, the function $ \gamma_t(x)$ also has a clipping procedure.  The true bandwidth is $10^{1/2}h_n/f^{1/2}(t)$ which is a constant if $f(X_i)< f(t)/10$.  Besides of the square root law, the clipping procedure is also necessary to improve the bias from the order of $h_n^2$ to the order of $h_n^4$ for the estimator (\ref{abr1}). The clipping procedures prevent too much contribution to the density estimation at $t$ if the observation $X_i$ is too far away from $t$. Later works on variable bandwidth density estimation include \cite{HM1988},  \cite{HallHuMarron},  \cite{TerrellScott}, \cite{McKay a}, \cite{McKay b},  \cite{Novak}, \cite{JonesMcKayHu}, \cite{GS2013, GineSang}, \cite{NS2016} and so on.

However,  this variable bandwidth estimator (\ref{abr1}) is not a density function since the integral of $f_n(t)$ over $t$ is not  $1$.  McKay \cite{McKay a, McKay b} discovered a smooth clipping function and studied the following variable bandwidth kernel density estimation 
\begin{align}\label{GSV}
\tilde{f}_n(t)= \frac{1}{n h_n}\sum_{i=1}^n K\left ( \frac{t-X_i}{h_n} \alpha \left (f(X_i) \right)  \right )  \alpha \left (f(X_i) \right). 
\end{align}
The smooth function in \cite{McKay a, McKay b} has form 
\begin{eqnarray}\label {afunction}
\alpha(w) = c p^{1/2}(w/c^2), 
\end{eqnarray}
where $c>0$ is a constant and the clipping function $p$ has at least fourth order derivative and satisfies $p(u) \geq 1$ for all $u$ and $p(u)=u$ for all $u \geq t_0$ for some $1 \leq t_0<\infty$. The function $\alpha(w)$ provides the square root law since $\alpha(w)=w^{1/2}$ if $w\ge t_0c^2$. It also provides the clipping procedure since $\alpha(w)\ge c$. This variable bandwidth estimator is a density function since the integral of $\tilde{f}_n(t)$ over $t$ is  $1$. See the study of this estimator in \cite{GS2013, NS2016}. Examples of clipping functions are given in \cite{McKay a} and \cite{GS2013}. 

Motivated by the work in variable kernel density estimation, in particular the idea of square root law and clipping procedure in the paper \cite{Abramson, McKay a, McKay b}, in this paper,  to improve the accuracy of Nadaraya-Watson estimator (\ref{hatrt}), we propose the following version of the variable bandwidth regression estimator,
\begin{eqnarray} \label {rnt}
\bar{r} (t; h_n) =  \frac{\bar{g}(t; h_n)}{\bar{f}(t;h_n)},
\end{eqnarray}
where
\begin{eqnarray} 
\bar{g}(t; h_n)= \frac{1}{n h_n}\sum_{i=1}^n K\left ( \frac{t-X_i}{h_n}  \alpha \left (q(X_i) \right)  \right ) \alpha \left (q(X_i) \right)    Y_i, \label {barg}\\
\bar{f}(t;h_n)= \frac{1}{n h_n}\sum_{i=1}^n K\left ( \frac{t-X_i}{h_n} \alpha \left (q(X_i) \right)  \right )  \alpha \left (q(X_i) \right)  . \label {barf}
\end{eqnarray}
Here
\begin{eqnarray}\label {VB55}
q(x) =f(x) \sqrt{|r'(x)|} .
\end{eqnarray} 
Notice that $\alpha(q(x))=q^{1/2}(x)$ if $q(x)\ge t_0c^2$ by the definition of $\alpha(\cdot)$ in (\ref{afunction}). Similar to the square root law in the variable density estimation in \cite{Abramson}, we use a variable bandwidth $h_n/q^{1/2}(X_i)$ at the observation $X_i$ to estimate $r(t)$ if $q(X_i)\ge t_0c^2$. The intuition is, if $f$ or $|r'|$ is large at observation $X_i$, because of the continuity of  $f(x)$ and $r'(x)$, one expects to have more observations in the neighborhood of that $X_i$, and one should choose a small bandwidth to prevent too much data involved in the estimation. If $f$ or $|r'|$ is small at observation $X_i$, one expects to have much less observations in the neighborhood of that $X_i$, and one should choose a relative large bandwidth to pick up more data in the estimation.  On the other hand, if $q$ is close to zero at obserbvation $X_i$, instead of the extremely large bandwidth $h_n/q^{1/2}(X_i)$, one should use a bandwidth proportional to $h_n$  to avoid over estimation. This is realized by the clipping procedure in $\alpha$ since $\alpha(q(x))\ge c$. As stated in our main results, this selection of the bandwidth $h_n/\alpha(q(X_i))$ at each observation $X_i$ results in a bias with the order of $h_n^4$. 

The estimator (\ref{rnt}) is called ideal estimator because the function $q(x)$ in (\ref{VB55}) involves the functions $f(x)$ and $r(x)$ which are to be estimated. To have a practical version, we shall take two sequences of bandwidth $h_{1,n}$ and $h_{2,n}$ with $h_{1,n}, h_{2,n}\rightarrow 0$ and $nh_{1,n}, nh_{2,n}\rightarrow \infty$ as $n\rightarrow \infty$. The first sequence $h_{1,n}$ is for the initial estimation of $q(x)$. Let 
\begin{eqnarray} \label{hatq}
\hat{q}(x; h_{1,n})= \hat{f}(x; h_{1,n}) \sqrt{\left |\hat{r}'(x; h_{1,n})\right |}
\end{eqnarray}
where $\hat{f}(x; h_{1,n})$ and $\hat{r}(x; h_{1,n})$ are defined as in (\ref{hatft}) and (\ref{hatrt}) with $h_n$ replaced by $h_{1,n}$ and $\hat{r}'(x; h_{1,n})=d \hat{r}(x; h_{1,n})/dx$. Define the true estimator 
\begin{eqnarray}  \label{hatrh1h2}
\hat{r}(t; h_{1,n}, h_{2,n})=\frac{\hat{g}(t; h_{1,n}, h_{2,n})}{\hat{f}(t; h_{1,n}, h_{2,n})},
\end{eqnarray}
where
\begin{eqnarray*} 
\hat{g}(t; h_{1,n}, h_{2,n})= \frac{1}{n h_{2,n}}\sum_{i=1}^n K\left ( \frac{t-X_i}{h_{2,n}} \alpha(\hat{q}(X_i; h_{1,n}))   \right )  \alpha(\hat{q}(X_i; h_{1,n}))    Y_i,  \label {VB46b} \\
\hat{f}(t;h_{1,n}, h_{2,n})= \frac{1}{n h_{2,n}}\sum_{i=1}^n K\left ( \frac{t-X_i}{h_{2,n}} \alpha(\hat{q}(X_i; h_{1,n}))   \right )  \alpha(\hat{q}(X_i; h_{1,n})) . \label {VB45b}
\end{eqnarray*}
The idea of variable bandwidth is to assign a large bandwidth in sparse area and a small bandwidth in dense area. In an extremely sparse area, one applies a bandwidth proportional to $h_n$ to avoid over estimation. This procedure is smooth due to the differentiability of the clipping function.

An incomplete list of study on variable bandwidth kernel regression estimation includes \cite{GM1984}, \cite{MS}, \cite{Novak}, \cite{Hall} and \cite{EM2005}.  Einmahl and Mason \cite{EM2005} worked on establishing consistence of kernel-type estimators $(\ref{rtq})$ in the multidimensional case when the bandwidth $h_n$ is a function of the location $t$ or the data. 
M\"uller and Stadtm\"uller \cite{MS} studied kernel regression estimation with fixed design points $\{x_i\}_{i=1}^n$ and variable bandwidth depending on the estimation point $t$.
Hall \cite{Hall} used the variable bandwidth $h_{n}/\alpha_1(X_i)$ for the numerator and $h_{n}/\alpha_2(X_i)$ for the denominator of the regression estimator where $\alpha_1=|fg|^{1/2}=f|r|^{1/2}$ and $\alpha_2=f^{1/2}$ were recommended. 

The paper is organized as follows. Section \ref{main} gives the main results.   We present a simulation study in Section \ref{simulation}.  Section \ref{conclusion} concludes the paper with a brief summary. All proofs are reserved to Section \ref{proof}.

\bigskip

\textbf{Acknowledgement} 
The authors are grateful to the referees and the Associate
Editor for carefully reading the paper and for insightful suggestions that significantly
improved the presentation of the paper.  
The research of Hailin Sang is supported in part by
the Simons Foundation Grant 586789 and the College of Liberal Arts
Faculty Grants for Research and Creative Achievement at the University of Mississippi. The research of Janet Nakarmi is supported by University Research Council (URC) Grant at the University of Central Arkansas.


\section {Main results}\label{main}


Let $D_k(\cdot)$ denote the $k^{th}$ order derivative for $k \geq 1$ and $D_0(f)=f$. For integers $k \geq 0$ and $p \geq 1$, denote
\begin{eqnarray*}
\mu_{k, p} = \int w^k K^p(w)dw.
\end{eqnarray*}
Denote $m(t)=E(Y_1^2|X_1=t)$ and $\sigma^2(t)=m(t)-r^2(t)$.  \\

We heavily rely on the following propositions in the proof of the main theorems in this section. These two propositions are modification and  generalization of the results in \cite{McKay b}, \cite{JonesMcKayHu}, \cite{GS2013} (see also \cite{Hall} or \cite{NS2016}). 
 The density function $f(s)$ of random variable $X$ there can be replaced by the function $\eta(s)$ with the same smoothing property. 

\begin{proposition} \label{bias0rx}
Suppose that $K$ has bounded support $[-T, T]$ and integrates to $1$. Assume that $\eta$ and $\xi$  have $l+1$ bounded and continuous derivative, $\xi\ge c>0$, and $\xi'/\xi$ is bounded in a neighborhood of $t$  for some $c>0$. Then
\begin{equation}\label{P1}
\frac{1 }{h_n}\int K\left ( \frac{t-s}{h_n} \xi(s) \right ) \xi(s) \eta(s) ds=\sum_{k=0}^l a_{k}(t)h_n^k+o(h_n^l),
\end{equation}
 as $h_n \to 0$,  where the set of functions $a_{k}(t)$ are defined as 
\begin{equation} \label {P2}
 a_{k}(t)=(-1)^k\frac{\mu_{k, 1}}{k!}D_{k}\left(\frac {\eta(t)}{\xi^{k}(t)}\right).
\end{equation}
If $K$ is symmetric with respect to zero, then $a_{2k+1}(t)=0$. 
\end{proposition}

\begin{proposition} \label{bias0} 
Under the same condition as in Proposition \ref{bias0rx},
\begin{equation*}\label{locationideal3}
\frac{1}{h_n}\int K^2\left(\frac{t-s}{h_n}\xi(s)\right)\xi^{2}(s)\eta(s)ds =\sum_{k=0}^l a_{k}(t)h_n^{k}+o(h_n^{l})
\end{equation*}
 as $h_n \to 0$. The functions $a_{k}(t)$ are defined as
\begin{equation*}\label{coeff}
a_{2k+1}(t)=0,\ \ a_{2k}(t)= \frac{\mu_{2k,2}}{(2k)!}D_{2k}\left( \frac{\eta(t)}{\xi^{2k-1}(t)}\right).
\end{equation*}
\end{proposition}

We take the following assumption on the kernel function $K$, the clipping function $p$, the density function $f(x)$ of $X_1$, the regression function $r$  throughout the paper.
\begin{assumption} \label{ass0}
Assume that $K$ is non-negative, symmetric to zero, integrates to $1$ and has support $[-T, T]$ for some $T<\infty$, $K$ has fourth order derivative, the clipping function $p$ (see the definition of function $\alpha$ in   (\ref{afunction})) has fifth order derivative, $f$ has fifth order bounded and continuous derivative in the neighborhood of $t$, $r$ has sixth order bounded and continuous derivative in the neighborhood of $t$. $f'/f$ and $r''/r'$ are bounded in the neighborhood of $t$.
\end{assumption}

\begin{remark}\label{xicon}
For the proof of the main results in this section, we shall apply Propositions \ref{bias0rx} and \ref{bias0} several times with $\xi(s)=\alpha(q(s))$ and $l=4$. Recall that $q(s)=f(s)|r'(s)|^{1/2}$. Obviously, by the definition of function $\alpha$ in (\ref{afunction}) and the condition on $p$, $f$ and $r$ in Assumption \ref{ass0},  $\xi\ge c$ and $\xi$ has fifth order bounded and continuous derivative. We now show that this $\xi$ also satisfies the condition that $\xi'/\xi$ is bounded in a neighborhood of $t$ under the condition $f'/f, r''/r'$ are bounded in a neighborhood of $t$. 
We only need to consider the case $r'(s)>0$.
 If $q(s) \geq t_0c^2$, then $\xi(s)=\alpha(q(s))=f^{1/2}(s) (r'(s))^{1/4}$. Hence
\begin{eqnarray*}
\frac{d\alpha(q(s))/ds}{\alpha(q(s))}
= \frac{f'(s)}{2f(s)}+ \frac{r''(s)}{4 r'(s)}.
\end{eqnarray*}
 If $q(s) < t_0c^2$, $\xi(s)=\alpha(q(s))=cp^{1/2}(q(s)/c^2)$, we have 
\begin{eqnarray*}
\frac{d\alpha(q(s))/ds}{\alpha(q(s))}=\frac{  q'(s)p'(q(s)/c^2)}{2c^2 p(q(s)/c^2)}.
\end{eqnarray*}
Here
\begin{eqnarray*}
q'(s)= (f(s)(r'(s))^{1/2})'= q(s) \left \{ \frac{f'(s)}{f(s)}+ \frac{r''(s)}{2 r'(s)} \right \}, 
\end{eqnarray*}
and $0\le q(s)<t_0c^2$ in this case. $p'(q(s)/c^2)$ is bounded in a neighborhood of $t$ because of the continuity of $p'$. Therefore for all $s$ in a neighborhood of $t$, the boundedness of  $\frac{f'(s)}{f(s)}$ and $\frac{r''(s)}{ r'(s)}$ implies that $\xi'/\xi$ is bounded in a neighborhood of $t$.
\end{remark}
\begin{remark}
In the proof of the theorems, we apply  Propositions \ref{bias0rx} and \ref{bias0} and take the function $\eta$ there to be $f$ or $g=fr$. It is easy to see that the conditions on $f$ and $r$ in Assumption \ref{ass0} imply that $\eta$ has $l+1$ bounded and continuous derivative with $l=4$. 
\end{remark}

\begin{remark}
Gin\'{e} and Sang  \cite{GS2013} provided a five time differentiable clipping function $p$ with $t_0=2$. 
\begin{displaymath}\label{p5}
p(t)=\left\{\begin{array}{ll}
1+\frac{t^6}{64}\left(1-2(t-2)+\frac{9}{4}(t-2)^2-\frac{7}{4}(t-2)^3+\frac{7}{8}(t-2)^4\right) & \textrm{if $0\le t\le2$}\\
t&\textrm{if $t\ge 2$}\\
1&\textrm {if $t\le 0$}\end{array}
\right..
\end{displaymath}
We will use this clipping function in the simulation study. 
\end{remark}
Let  $1 \leq t_0 < \infty$ and $0<c<\infty$ be the constants in the definition of $\alpha(w)$ in (\ref{afunction}), we study the estimation of $r(t)$ for $t$ in the region ${\cal D}_{rf}$, 
\begin{eqnarray}\label{Drf}
{\cal D}_{rf}= \{ t \in \mathbb{R}: q(t) \geq   2 t_0 c^2    \}.
\end{eqnarray}
Note that $f(t)$ is bounded away from zero for $t \in D_{rf}$. This is a necessary
condition for the estimation of $r(t)$.
\begin{remark}
By the condition $q(t)=f(t)|r'(t)|^{1/2}\geq   2 t_0 c^2$ in the definition of region ${\cal D}_{rf}$, $\alpha(q(t))=q^{1/2}(t)$. This is a necessary requirement to remove the term with $h_n^2$ in the bias expansion. For example, see the proof of the following Theorem \ref{Thm01}, Step 1. On the other hand, we can also observe this point from the simulation study in Section \ref{simulation}. The estimator does not have satisfied performance in the area where $r'$ is close to $0$ (the area with flat tangent lines). See Figures \ref{fig2} and \ref{fig3} in Section \ref{simulation}.
\end{remark}

Now we are ready to state the main theorems. The next two theorems are the results for the ideal estimator in (\ref{rnt}). 
\begin{theorem}\label {Thm01}
Under Assumption \ref{ass0}, assume that $h_n\to 0$ and $nh_n\to \infty$ as $n\to \infty$, then 
\begin{eqnarray}\label {Thm1(1)}
E (\bar{r}(t; h_n)) = r(t) +  \theta(t) h_n^4   +o(h_n^4)+o\left (\frac{1}{n h_n} \right ) 
\end{eqnarray}
and 
\begin{eqnarray}\label {Thm1(2)}
E\left ( \bar{r}(t; h_n)- r(t) \right )^2 =\theta^2(t) h_n^8 +\frac{\mu_{0,2}|r'(t)|^{1/4}\sigma^2(t)}{ nh_n \sqrt{f(t)}}  +o(h_n^8)+o\left (\frac{1}{n h_n}  \right)
\end{eqnarray}
for $t \in {\cal D}_{rf}$, where
\begin{eqnarray}\label {VB12a}
\theta(t) = \frac{\mu_{4,1}   }{24  f(t)}\left \{ D_4 \left ( \frac{r(t) }{f(t) |r'(t)|} \right ) - r(t) D_4 \left ( \frac{ 1}{f(t)|r'(t)|} \right )\right \} .
\end{eqnarray}
\end{theorem}


\begin{theorem}\label {central1}
Under Assumption \ref{ass0},  if $h_n^4 \sqrt{n h_n} \rightarrow \lambda$ as $n\to\infty$, for some $0 \leq \lambda <\infty$, then
\begin{eqnarray*}
\sqrt{n h_n} \left \{\bar{r}(t; h_n)-r(t) \right\} \xrightarrow{D} N \left (\lambda \theta(t),  \frac{\mu_{0,2}|r'(t)|^{1/4}\sigma^2(t)}{\sqrt{f(t)} }  \right )
\end{eqnarray*}
for $t \in {\cal D}_{rf}$.
\end{theorem}

\begin{remark}
M\"uller and Stadtm\"uller \cite{MS} worked on kernel regression estimation with fixed design points $\{x_i\}_{i=1}^n$ and variable bandwidth depending on the estimation point $t$. In this paper we study kernel regression estimation with random design points $\{X_i\}_{i=1}^n$ and variable bandwidth depending on the sample. 

Hall \cite{Hall} used the variable bandwidth $h_{n}/\alpha_1(X_i)$ for the numerator and $h_{n}/\alpha_2(X_i)$ for the denominator of the regression estimator where $\alpha_1=|fg|^{1/2}=f|r|^{1/2}$ and $\alpha_2=f^{1/2}$. However, if we write Hall's regression estimator as $\sum_{i=1}^n w_iY_i$, the sum of the weight $\sum_{i=1}^n w_i\ne 1$ since the bandwidths for the numerator and the denominator are different. In this paper we use the same variable bandwidth $h_n/\alpha(q(X_i))$ for both the numerator and the denominator. Consequently, if we define the weight $w_i=(n h_n)^{-1}K\left ( \frac{t-X_i}{h_n}  \alpha \left (q(X_i) \right)  \right ) \alpha \left (q(X_i) \right)/\bar{f}(t; h_n)$ for each $X_i, i \in[1,n]$, then $\sum_{i=1}^n w_i=1$.

On the other hand, Hall's regression estimator has no clipping procedure. As Terrell and Scott \cite{TerrellScott} pointed out in variable kernel density estimation, one may not have bias reduction if we only use the square root law. This may apply to variable bandwidth kernel regression estimation if the estimation of the density function is involved. 
\end{remark}

The next two theorems are the results for the true estimator in (\ref{hatrh1h2}). 
\begin{theorem} \label {Thm2}
Denote $U(h_{1,n}):=h_{1,n}^2+( n h^3_{1,n})^{-1/2}$ and assume  $h_{2,n}\to 0$ and $nh_{2,n}\to \infty$ as $n\to \infty$, $U(h_{1,n})=o(h_{2,n}^2)$. Under Assumption \ref{ass0}, 
\begin{eqnarray} 
E \hat{r}(t; h_{1,n}, h_{2,n}) =r(t) +\theta(t) h^4_{2,n}+o(h^4_{2,n})+o\left ( \frac{1}{n h_{2,n}} \right)  \label {Thm2(1)} 
\end{eqnarray} 
and 
\begin{eqnarray} 
&&E \left (\hat{r}(t; h_{1,n}, h_{2,n})-r(t) \right)^2 \nonumber\\
&&= \theta^2(t) h_{2,n}^8 +\frac{|r'(t)|^{1/4}\mu_{0,2}\sigma^2(t)}{ nh_{2,n}\sqrt{f(t)}}   + o(h^8_{2,n})+o\left (\frac{1}{n h_{2,n}}  \right )\label {Thm2(2)}
\end{eqnarray} 
for $t \in {\cal D}_{rf}$. Consequently, the optimal bandwidth 
\begin{eqnarray*}
h_{2,n}^*=\left (\frac{1}{n}  \right )^{1/9}\left ( \frac{\mu_{0,2}\sigma^2(t)}{8 \int_{t \in D_{rf}}\theta^2(t)dt}\int_{t \in D_{rf}}\frac{|r'(t)|^{1/4}}{\sqrt{f(t)}} dt\right )^{1/9}.
\end{eqnarray*}
\end{theorem}


\begin{theorem}\label {central2}
Assume $U(h_{1,n})=o(h_{2,n}^2)$. Under Assumption \ref{ass0}, if $h_{2,n}^4 \sqrt{n h_{2,n}} \rightarrow \lambda$ as $n\to\infty$,  for some  $0 \leq \lambda <\infty$, 
\begin{eqnarray*}
\sqrt{n h_{2,n}} \left \{\hat{r}(t; h_{1,n}, h_{2,n}) -r(t) \right\} \xrightarrow{D} N \left (\lambda \theta(t),  \frac{\mu_{0,2}|r'(t)|^{1/4}\sigma^2(t)}{\sqrt{f(t)} }  \right )
\end{eqnarray*}
for $t \in {\cal D}_{rf}$. 
\end{theorem}


\section {Simulation}\label{simulation}
\indent In this section, we conduct simulation study to compare the performance of the variable bandwidth kernel regression estimator (VKRE) (\ref{hatrh1h2}) with Nadaraya Watson estimator (NWE). The three regression functions used are:
\begin{eqnarray} \label{reg_fun}
&1.&\;\;Y_i = 2+ \sin(0.75X_i)+0.3 \varepsilon_i; \;\;\;2.\;\;Y_i =\frac{1}{1+X_i^2}+0.3 \varepsilon_i ; \nonumber\\
&3.&\;\;Y_i =\log{|X_i|}+0.3 \varepsilon_i, 
\end{eqnarray}

\noindent where $1\le i\le n$, $\{\varepsilon_i\}_{i=1}^n $ are $i.i.d.$ random errors, $\{X_i\}_{i=1}^n$ are $i.i.d.$ random variables, and they are independent.  For NWE, we used the  {\it npreg()} function from  the {\it np} package \cite{Hay} in the programming software R with their default settings, i.e., Gaussian kernel and cross-validation bandwidth selection method. For VKRE, we applied the following kernel (tricube kernel)  
\begin{equation*}
K(u) = \frac{70}{81}(1-|u|^3)^3 1_{|u| \le 1},
\end{equation*}
and the five time differentiable clipping function $p$ in (\ref{p5}). 

Recall that we estimate the regression function $r(t)$ for $t \in {\cal D}_{rf}$.   Hence, we can estimate $r(t)$ over a large range of $t$ if we choose a very small value of $c$. In this simulation study we set  $c$ in (\ref{afunction})  to be $0.000001$. To apply the kernel method in an estimation, one should select an optimal bandwidth based on some criteria, for example, to minimize the mean squared error. We apply the cross-validation bandwidth selection method in NWE.  It is interesting to investigate the bandwidth selection problem from both theoretical and application viewpoints for VKRE. However, the study in this direction is a new challenge and we leave it as an open question for future study. Instead, we take the bandwidths, $h_{1,n}$ and $h_{2,n}$,  as $0.6\times n^{-1/7}$ and $n^{-1/9}/4$ respectively for all following simulation study including the graphs and numerical comparison, which satisfy the assumptions in Theorem \ref{Thm2}. 

\begin{figure}[H]
\begin{center}
\includegraphics[ height=1.8in, width=4in]{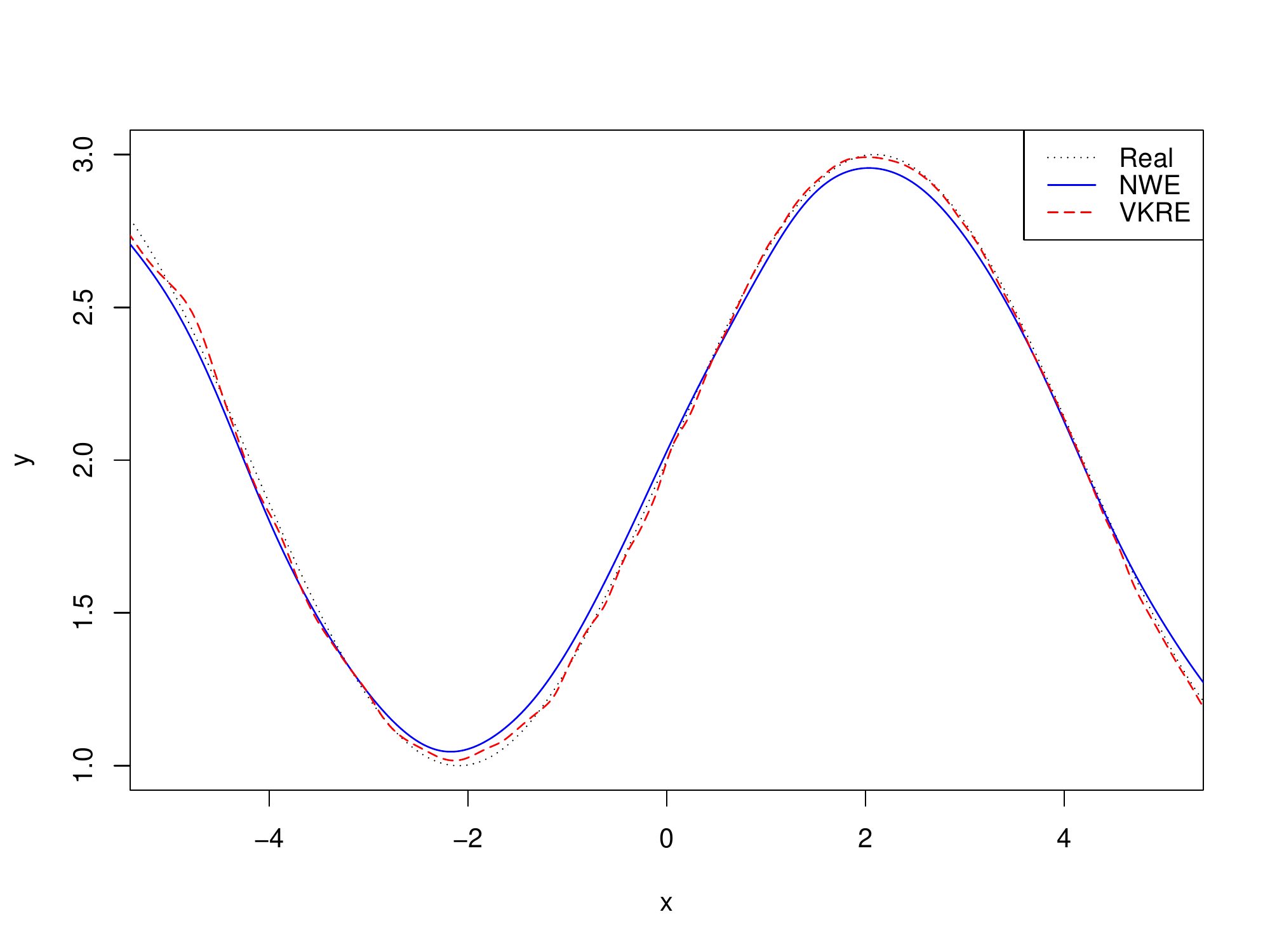}
\end{center}	 
	\caption{$y=2+ \sin(0.75x)$}
\label{fig2}
\end{figure}

\begin{figure}[H]
\begin{center}
\includegraphics[ height=1.8in, width=4in]{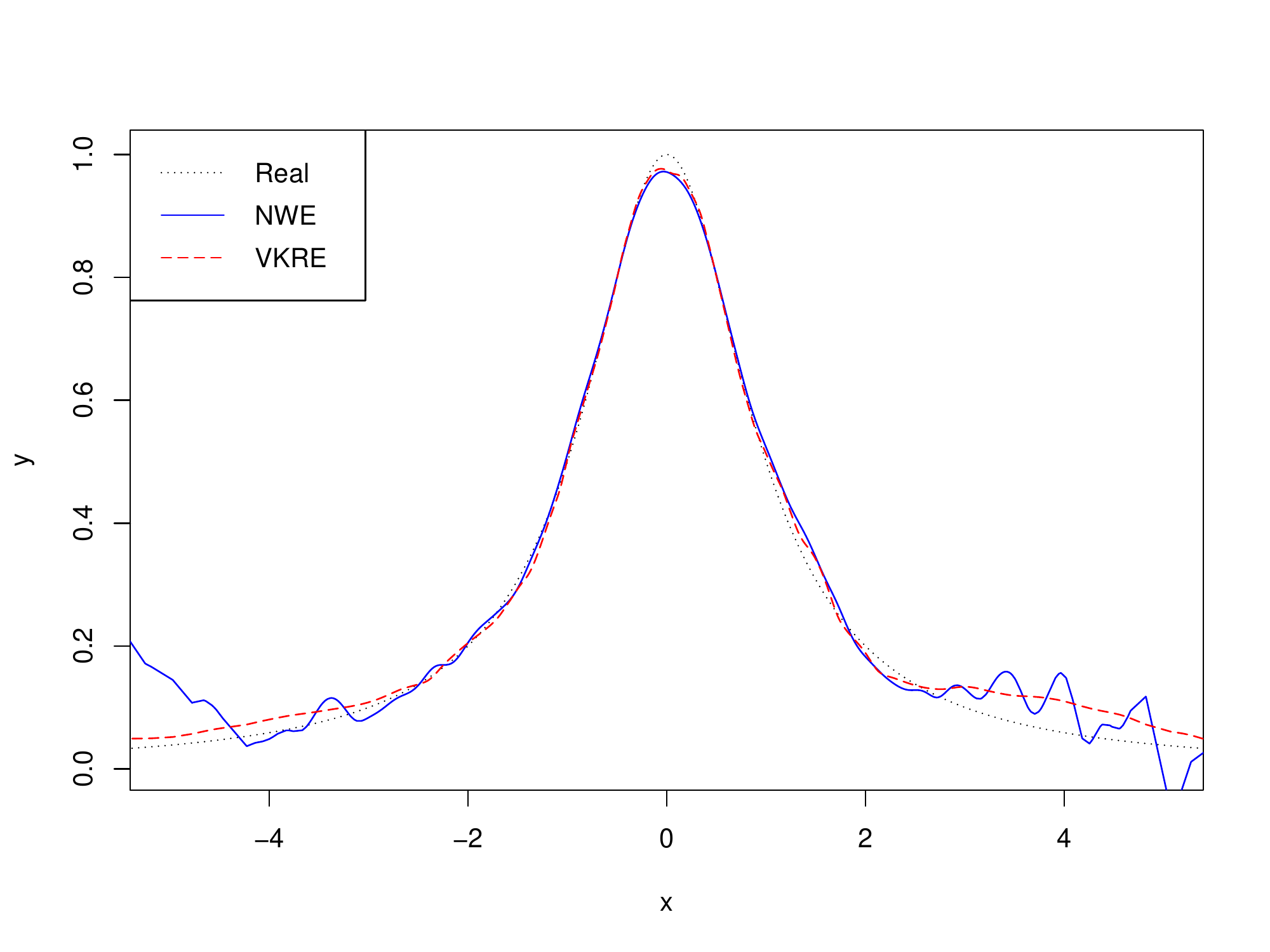}
\end{center}	 
	\caption{$y=1/(x^2+1)$}
\label{fig3}
\end{figure}

\begin{figure}[H]
\begin{center}
\includegraphics[ height=1.8in, width=4in]{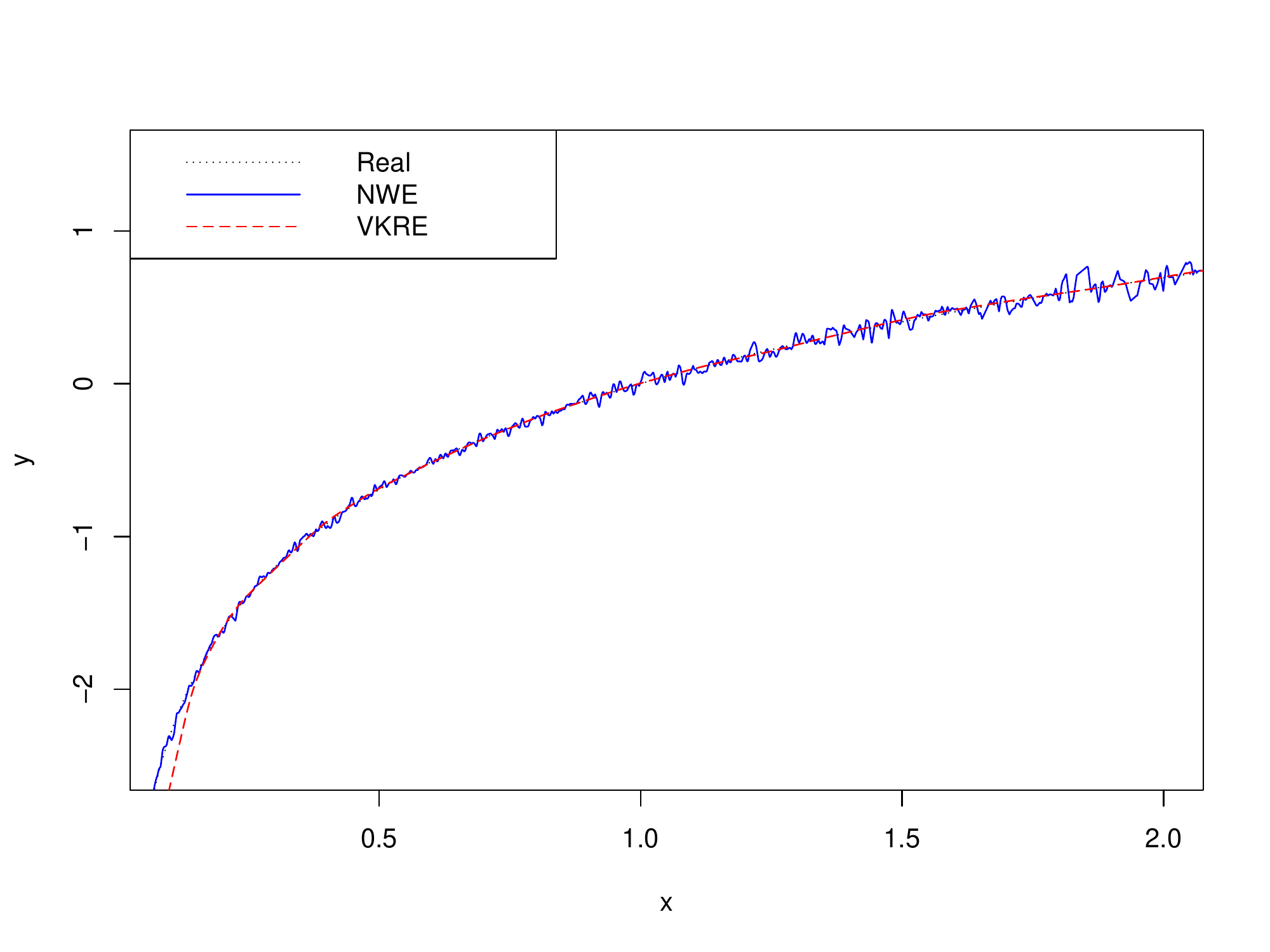}
\end{center}	 
	              \caption{$y=\log|x|$}
\label{fig4}
\end{figure}

For the simulation study in Figures \ref{fig2}, \ref{fig3}, and  \ref{fig4}, the random errors are generated from the uniform distribution on interval $[-0.5,0.5]$ with $n=5000$.   In Figure \ref{fig2}, the first regression function in (\ref{reg_fun}) is chosen with the random variables $\{X_i\}_{i=1}^n$  generated from the Cauchy distribution with location parameter $3$ and scale parameter $4$.  In Figure \ref{fig3}, the second regression function in (\ref{reg_fun}) is chosen with $\{X_i\}_{i=1}^n$ generated from the T-distribution with degree of freedom $4$. In Figure \ref{fig4}, the third regression function in (\ref{reg_fun}) is chosen with $\{X_i\}_{i=1}^n$ generated from the standard normal distribution. 

The simulation study in Figures \ref{fig2}, \ref{fig3}, and  \ref{fig4}  shows that, for each of the regression functions,  VKRE has better performance than NWE. The only exception is the area around $0$ in the third figure. For the regression function $r(x)=\log |x|$, if $x>0$ is in a neighborhood of $0$, $r''/r'=-1/x$ is not bounded. Hence the condition on $r$ in Assumption \ref{ass0} is not satisfied. This is the reason VKRE has very bad performance around $0$.   Also we notice that the performance of VKRE is slightly worse in the area $r'$ is close to $0$ ( the area with horizontal tangent line) than the other  part in Figures \ref{fig2} and \ref{fig3}. In some sense this confirms the condition $q(t)=f(t)|r'(t)|^{1/2}\ge 2t_0c^2$ in the definition of the region ${\cal D}_{rf}$ in (\ref{Drf}).  

 \subsection{Numerical Comparison using RMSE}
We provide numerical comparison to further verify the advantage of VKRE. Tables \ref{tab1}, \ref{tab2},  and \ref{tab3}  demonstrate the difference  between NWE and VKRE for the second regression function in (\ref{reg_fun}) using the measure,
$$RMSE= \sqrt{\frac{1}{n}\sum_{i=1}^n \big(r(X_i) - \hat{r}(X_i)\big)^2},$$
 where $r(t)$ is the real regression function and $\hat{r}(t)$ is NWE or VKRE. Each entry of the tables is the average of RMSE from $N=250$ samples. 
 
In Table \ref{tab1}, we compare the RMSE for NWE and VKRE with the random errors from different distributions, $X$ from the T-distribution with degree of freedom $4$, and the sample size of $5000$. The table shows that as the bounds for the uniform distribution increases, the RMSE increases for both NWE and VKRE.

 \begin{table}[H]
\center
\caption{ Comparing  RMSE of VKRE and NWE with the random errors from different distributions.} 
\label{tab1}
\bigskip
\begin{tabular}{| c| c| c| }
\hline                                                     
\textbf{\;\;Errors\;\;}   &  \textbf{\;\;NWE RMSE\;\;} &  \textbf{\;\;VKRE RMSE\;\;} \\
\hline
$U[-0.5,0.5]$ & 0.01165791& 0.008649485\\
\hline
$U[-1,1]$ & 0.02135706& 0.01474327\\
\hline
$U[-2,2]$ & 0.03912875& 0.02785478\\
\hline 
\end{tabular}  
\end{table}

  In Table \ref{tab2}, we list the RMSE for NWE and VKRE with $X$ generated from different distributions, $\varepsilon$ generated from the uniform distribution on the interval $[-0.5,0.5]$, and $n=5000$ with $N=250$ repetitions. Recall that the T-distribution with degree of freedom $1$ is same as the standard Cauchy distribution. It is interesting to see that as the degrees of freedom increase for T-distribution, the RMSE decreases for both of the estimators.

  \begin{table}[H]
\center
\caption{ Comparing  RMSE of VKRE and NWE with $X$ from different distributions.} 
\label{tab2}
\bigskip
\begin{tabular}{| c| c| c| }
\hline                                                     
\textbf{\;\;$X$\;\;}   &  \textbf{\;\;NWE RMSE\;\;} &  \textbf{\;\;VKRE RMSE\;\;} \\
\hline
$T(df=1)$ & 0.01973383& 0.01215381\\
\hline
$T(df=4)$ & 0.01165791& 0.008649485\\
\hline
$T(df=8)$ & 0.006720047 & 0.004981605 \\
\hline
$Cauchy (3,4)$ & 0.02445313 &0.01189971 \\
\hline 
$Cauchy (5,7)$ & 0.02656666 & 0.01461736 \\
\hline 
$N(0,1)$ & 0.007069801 &0.005971244 \\
\hline 
$N(5,10)$ & 0.01349415 & 0.01335264 \\
\hline 
\end{tabular}  
\end{table}

 In Table \ref{tab3}, we compare the RMSE for NWE and VKRE with different sample sizes, $\varepsilon$ generated from the uniform distribution on the interval $[-0.5,0.5]$, and $X$ generated from the T-distribution with degree of freedom $4$. The table shows that the RMSE decreases as the sample size increases for both of the estimators.

\begin{table}[H]
\center
\caption{ Comparing RMSE of VKRE and NWE with different sample sizes.} 
\label{tab3}
\bigskip
\begin{tabular}{| c| c| c| }
\hline                                                     
\textbf{\;\;$n$\;\;}   &  \textbf{\;\;NWE RMSE\;\;} &  \textbf{\;\;VKRE RMSE\;\;} \\
\hline
$500$ & 0.01866795 & 0.01520967\\
\hline
$1000$ & 0.01598541 & 0.01162556 \\
\hline
$2000$ & 0.01358658 &  0.01012766\\
\hline
$5000$ & 0.01165791& 0.008649485\\
\hline 
$8000$ &  0.009052834 & 0.00564806  \\
\hline 
$10000$ & 0.007793899 &0.004173846  \\
\hline 
\end{tabular}  
\end{table}
Moreover, the results from all three tables indicate that  the RMSE of VKRE is smaller than that of NWE in each case. Thus, in all situations considered,  VKRE outperforms  NWE.

 \subsection{Numerical Comparison using Monte-Carlo Estimation of MSE}
 This subsection covers Monte-Carlo simulation and approximation of the following mean square error,
 $$MSE= \mathbb{E}[ \big(\hat{r}(t) - r(t)\big)^2],$$
 to compare NWE and VKRE estimators.  We compare the regression estimators for the second regression function (Bounded RF) and the third regression function (Unbounded RF) in (\ref{reg_fun}) in the following tables. Each entry in the tables is the average of $(\hat{r}(t) - r(t))^2$ from $N=250$ random samples each with sample size $n= 5000$. In each table, we list the results at $10$ points which are evenly selected from the range of $nN$ simulated $X$ values with the specified T-distribution or the standard normal distribution. \\
 
Tables \ref{tab4} and \ref{tab5} compares NWE and VKRE for the bounded and the unbounded regression functions in (\ref{reg_fun})  where $X$ is generated from the T-distribution with degrees of freedom 4, whereas in Tables \ref{tab6} and \ref{tab7},  $X$ is generated from the standard normal distribution.   In Tables \ref{tab4} and  \ref{tab6} the errors are from the standard normal distribution whereas the errors in Tables  \ref{tab5} and  \ref{tab7} are generated from the uniform distribution on the interval [-1,1].

 \begin{table}[H]
	\center
	\caption{ Comparing  MSE of VKRE and NWE with $X \sim T(4)$ and $\varepsilon \sim N(0,1)$.} 
	\label{tab4}
	\bigskip
\begin{tabular}{|c|l|l|l|l|l|l|}
\hline
 &
\multicolumn{2}{c}{ \textbf{Bounded RF}} &
\multicolumn{2}{c|}{ \textbf{Unbounded RF}}  \\
\hline
\textbf{\;\;$t$\;\;}   &  \textbf{\;\;NWE \;\;} &  \textbf{\;\;VKRE \;\;} &  \textbf{\;\;NWE \;\;} &  \textbf{\;\;VKRE \;\;}\\
\hline
-7.161518 & 0.082425068 & 0.022757438 & 0.952882068  &   0.012199782\\
\hline
-5.593896  & 0.038621453 &0.029278125 &   0.181799310 & 0.010549010  \\
\hline
  -4.026274 & 0.001274921 &0.001691833 & 0.059748939 &  0.000630284 \\
\hline
 -2.458652& 0.000693960 & 0.000953009 &  0.028731527 & 0.000743185  \\
\hline
 -0.89103& 0.000095033 & 0.000350961 & 0.003600552 & 0.000450090  \\
\hline
 0.676592& 0.000470135 &0.000042570 &  0.000808568 & 0.000039578 \\
\hline
 2.244214 & 0.006266251 & 0.002081689 &  0.027084148 & 0.002896018 \\
\hline
 3.811836 & 0.006611614 & 0.000129696 &  0.049252162 &  0.002134087 \\
\hline
 5.379458 &  0.021020283 & 0.001573274 &  0.067687352 &  0.014088619 \\
\hline
 6.94708 &0.013454439  &  0.004420131&  0.709749776 &  0.042635434 \\
\hline
\end{tabular}
\end{table}

\begin{table}[H]
	\center
	\caption{ Comparing  MSE of VKRE and NWE with $X \sim T(4)$ and $\varepsilon \sim U[-1,1]$.} 
	\label{tab5}
	\bigskip
\begin{tabular}{|c|l|l|l|l|l|l|}
\hline
 &
\multicolumn{2}{c}{ \textbf{Bounded RF}} &
\multicolumn{2}{c|}{ \textbf{Unbounded RF}}  \\
\hline
\textbf{\;\;$t$\;\;}   &  \textbf{\;\;NWE \;\;} &  \textbf{\;\;VKRE \;\;} &  \textbf{\;\;NWE \;\;} &  \textbf{\;\;VKRE \;\;}\\
\hline
 -7.161518 &  0.037688485 &  0.015996216 & 1.4541232  &  0.0888033  \\
\hline
-5.593896 & 0.017563571  & 0.003636569 &   0.2238683 & 0.0268285  \\
 \hline
 -4.026274 & 0.000725586 & 0.000042321 &  0.0240457  &   0.0029750 \\
 \hline
 -2.458652 & 0.000536810 & 0.000052783 & 0.0103590  & 0.0001132  \\
 \hline
 -0.89103 & 0.000083700 & 0.000015774 & 0.0020801  &   0.0000299\\
 \hline
 0.676592 &  0.000015691 &  0.000047977 & 0.0019164   & 0.0000508  \\
 \hline
 2.244214 &  0.000151535 &  0.000032906 &   0.0111435 &  0.0002437 \\
 \hline
 3.811836 & 0.004435139 &  0.000604593 & 0.0250448  &  0.0005697 \\
 \hline
 5.379458 & 0.016278344  & 0.004865297 & 0.1074176  &  0.0271930  \\
 \hline
 6.94708 &  0.024189858 & 0.011530802  &   1.4184324 &   0.0756486\\
  \hline
\end{tabular}
\end{table}

In Tables  \ref{tab4}-\ref{tab7}, VKRE outperforms NWE for the unbounded regression function at all selected $t$ values.   For the bounded regression function,  VKRE has better performance than NWE at most of the $t$ values. 
On the other hand, for the $t$ values where NWE has better performance than VKRE, the difference between the corresponding MSEs is relatively very small. In summary, the  Monte-Carlo approximation of MSE also shows that VKRE outperforms NWE in general.

 \begin{table}[H]
	\center
	\caption{ Comparing  MSE of VKRE and NWE with $X \sim N(0,1)$ and $\varepsilon \sim N(0,1)$.} 
	\label{tab6}
	\bigskip
\begin{tabular}{|c|l|l|l|l|l|l|}
\hline
 &
\multicolumn{2}{c}{ \textbf{Bounded RF}} &
\multicolumn{2}{c|}{ \textbf{Unbounded RF}}  \\
\hline
\textbf{\;\;$t$\;\;}   &  \textbf{\;\;NWE \;\;} &  \textbf{\;\;VKRE \;\;} &  \textbf{\;\;NWE \;\;} &  \textbf{\;\;VKRE \;\;}\\
\hline
 -3.166296000 &  0.044502685 & 0.014599755 &   0.162391619 &  0.001478955 \\
\hline
 -2.476748778 & 0.002788421 & 0.001986454 & 0.052629439  & 0.000617446  \\
 \hline
 -1.787201556 &  0.000879162 & 0.000831587 &  0.008172456  &  0.000803430 \\
 \hline
 -1.097654333 & 0.000227731  & 0.000011475  &  0.001634166 & 0.000075463  \\
 \hline
 -0.408107111 & 0.000120885 &  0.000182900 &  0.001949065  &   0.000539604 \\
 \hline
 0.281440111 & 0.000061373 &  0.000206260 &  0.000934502 & 0.000086808  \\
 \hline
 0.970987333 &  0.000007520 &  0.000083182 & 0.001005907   &  0.000117938 \\
 \hline
 1.660534556 & 0.003158112 &  0.001136445 &   0.017263192 &  0.000629518  \\
 \hline
 2.350081778 &  0.000985980 &  0.000025147 & 0.035076566  &  0.003724841 \\
 \hline
 3.039629000 & 0.003365033 &  0.000184329 & 0.026700199  &  0.025010083 \\
  \hline
\end{tabular}
\end{table}

\begin{table}[H]
	\center
	\caption{ Comparing  MSE of VKRE and NWE with $X \sim N(0,1)$ and $\varepsilon \sim U[-1,1]$.} 
	\label{tab7}
	\bigskip
\begin{tabular}{|c|l|l|l|l|l|l|}
\hline
 &
\multicolumn{2}{c}{ \textbf{Bounded RF}} &
\multicolumn{2}{c|}{ \textbf{Unbounded RF}}  \\
\hline
\textbf{\;\;$t$\;\;}   &  \textbf{\;\;NWE \;\;} &  \textbf{\;\;VKRE \;\;} &  \textbf{\;\;NWE \;\;} &  \textbf{\;\;VKRE \;\;}\\
\hline
 -3.166296000 &  0.030718309 &  0.002079258 & 0.067146594  &  0.051639374 \\
\hline
 -2.476748778 &  0.000201962 &  0.000035443 &   0.018435422 &  0.006108237 \\
 \hline
 -1.787201556 &  0.000143322 &  0.000067562 &  0.003757715  &  0.000101261  \\
 \hline
 -1.097654333 & 0.000091673 & 0.000189295 &   0.002869324 & 0.000185627   \\
 \hline
 -0.408107111 & 0.000030893 &  0.000017899 & 0.001339005  & 0.000052513  \\
 \hline
 0.281440111 & 0.000013158  & 0.000017005 &  0.001477281  &  0.000069223  \\
 \hline
 0.970987333 & 0.000205468 & 0.000050455 & 0.001344766  & 0.000136607   \\
 \hline
 1.660534556 & 0.000036136 & 0.000004868 &  0.004016354  & 0.000039407  \\
 \hline
 2.350081778 &  0.001747355 &  0.000059404 & 0.019114022  & 0.001437631  \\
 \hline
 3.039629000 & 0.013097496 &  0.000167359 & 0.041389033  & 0.035245030   \\
  \hline
\end{tabular}
\end{table}


\section{Conclusion}\label{conclusion}

In this article we propose a variable bandwidth kernel regression estimator. With this estimator, the bandwidth is proportional to $1/\sqrt{f(X_i)}$, $1\le i\le n$, the inverse of the square root of the marginal density function value of the independent variable at the observation which is same as the square root law in the variable bandwidth kernel density estimation (\cite{Abramson}, \cite{HM1988},  \cite{HallHuMarron},  \cite{TerrellScott}, \cite{McKay a}, \cite{McKay b}, \cite{JonesMcKayHu}, \cite{GineSang}, \cite{GS2013} and \cite{NS2016}). The bandwidth is also proportional to the inverse of the absolute value of the fourth root of the derivative of the regression function.  It intuitively sounds since there are much more observations in the area with a large marginal density of the independent variable or with a large derivative of the regression function and we therefore shall take a relatively small bandwidth.  On the other hand, this variable bandwidth method also selects a larger bandwidth in the sparse area to pick up more observations in the estimation. In the area where the marginal density of the independent variable is extremely small or the regression is very flat, the clipping procedure there will take a relative constant bandwidth to avoid over estimation since the the clipping function is bounded below. 

Under some regular conditions on the kernel function and bandwidth sequence, we study the bias and mean squared error for both the ideal estimator and the true estimator. The order of the bias is $h_n^4$ instead of $h_n^2$ as in the classical Nadaraya-Watson kernel regression estimator. In consequence, the mean squared error has order of $n^{-8/9}$ instead of $n^{-4/5}$. We also obtain central limit theorem for the ideal estimator and the true estimator. The advantage of this variable bandwidth estimator over the classical Nadaraya-Watson kernel regression estimator is demonstrated by a simulation study.

It is also interesting to study the case when $X_i$'s are random variables on $\mathbb{R}^d$. In \cite{GS2013}, the authors studied variable bandwidth kernel density estimation (\ref{GSV}) in $d$-dimensional case. However, the true bandwidth in multidimensional variable bandwidth regression estimation should also involve the regression function $r$. It is an interesting problem to find the right variable bandwidth which can remove the term with $h_n^2$. We leave this part of work for future research.

\section {Proofs}\label{proof}



For convenience,  we slightly  modify the proofs of the theorems in \cite{McKay b}, \cite{JonesMcKayHu} or \cite{GS2013} and provide a proof for Proposition  \ref{bias0rx}.  The proof of Proposition  \ref{bias0} is similar. \\

\noindent\begin{proof}[Proof of Proposition  \ref{bias0rx}]
Since $\xi\geq c$ and  $\xi'/\xi$ is bounded in a neighborhood of $t$, there exists $\delta>0$ such that $(v\xi(t-v))'=\xi(t-v)-v \xi'(t-v)>0$ for $v \in [-\delta, \delta]$. Hence the function $U_t(v):= v \xi(t-v)$ is invertible for $v \in [-\delta, \delta]$. The inverse function $V_t(u)$ is $l+1$ times differentiable with continuous derivatives. Since $K((t-s)\xi(s)/h_n)=0$ unless $|t-s|\leq T h_n/\xi(s) \leq T h_n/c$, the change of variables 
\begin{eqnarray*}\label {hz}
h_n z = (t-s)\xi(t-(t-s)), \ \text{or} \ t-s = V_t(h_nz)
\end{eqnarray*}
in the following integral is valid 
\begin{eqnarray*}\label {Vhz}
&& \frac{1 }{h_n}\int K\left ( \frac{t-s}{h_n} \xi(s) \right ) \xi(s) \eta(s) ds \nonumber\\
&=&-\int K(z) \xi \left (t-V_t(h_n z)\right)\eta \left (t-V_t(h_nz)\right)\frac{dV_t(h_nz)}{d (h_nz)}dz. 
\end{eqnarray*}
Developing $\xi \left (t-V_t(h_n z)\right)\eta \left (t-V_t(h_nz)\right)\frac{dV_t(h_nz)}{d (h_nz)}$  into powers of $h_nz$, and the first statement of the proposition follows. 

Let $\psi$ be an infinitely differentiable function of finite support. Changing variable $t=s+h_n u$, developing $\psi$, changing variable $\omega = u \xi(s)$, and integrating by parts, we have
\begin{eqnarray}\label {23aa}
&& \frac{1 }{h_n}\int \psi(t) K\left ( \frac{t-s}{h_n} \xi(s) \right ) \xi(s) \eta(s) ds \nonumber\\
&=& \int \psi(s)\eta(s)ds + \sum_{k=1}^l (-1)^k \frac{\tau_k h_n^k}{k!}\int \psi(s) D_k \left (\frac{\eta(s)}{\xi(s)}  \right ) ds +o(h_n^l)
\end{eqnarray}
where $\tau_k=0$ when $k$ is odd by symmetry. By (\ref{P1}),
\begin{eqnarray}\label {24aa}
 \frac{1 }{h_n}\int \psi(t) K\left ( \frac{t-s}{h_n} \xi(s) \right ) \xi(s) \eta(s) ds =\sum_{k=0}^l h_n^k \int \psi(t) a_k(t)dt +o(h_n^l).
\end{eqnarray}
Comparing the coefficients of (\ref{23aa}) and (\ref{24aa}), we obtain (\ref{P2}). 
\end{proof}




We shall use the following formula to estimate the expectation of a quotient of two random variables:
\begin{eqnarray}\label {VB16ag}
\frac{1}{z}= 1-(z-1)+\cdots +(-1)^p (z-1)^p +(-1)^{p+1} \frac{(z-1)^{p+1}}{z}.
\end{eqnarray}


\noindent\begin{proof}[Proof of Theorem \ref{Thm01}]

By (\ref{rnt}) and (\ref{VB16ag}) with $p=1$ and $z=\bar{f}(t;h_n)/E \bar{f}(t;h_n)$,
\begin{eqnarray}   \label {21ta}
E \bar{r} (t; h_n)= \frac{E \bar{g}(t; h_n)}{E \bar{f}(t; h_n)} \label {VB32a}  +\frac{-I_1+I_2}{(E \bar{f} (t; h_n))^2} 
\end{eqnarray}
where
\begin{eqnarray}
I_1 &=& E\left \{\bar{g}(t;h_n) \left (\bar{f}(t;h_n)-E \bar{f}(t;h_n) \right) \right\}, \label {I_1}\\
I_2 &=& E \left \{\bar{r}(t;h_n) \left (\bar{f}(t;h_n)-E \bar{f}(t;h_n)\right)^2 \right\}.\label {I_2}
\end{eqnarray}
\textbf{Step 1}. We estimate $E \bar{g}(t;h_n)/E\bar{f}(t;h_n)$. Let $g(t)=f(t)r(t)$. Then
\begin{eqnarray}\label {29L}
\frac{E \bar{g}(t; h_n)}{E \bar{f}(t; h_n)}-r(t)= \frac{f(t) E \bar{g}(t; h_n)-g(t)E \bar{f}(t; h_n)}{f(t)E \bar{f}(t; h_n)}.
\end{eqnarray}
Note that
\begin{eqnarray*}\label {28wv}
E \bar{g}(t; h_n)= \frac{1}{h_n} \int K \left ( \frac{t-s}{h_n} \alpha \left (q(s) \right)  \right )\alpha(q(s))g(s)  ds.
\end{eqnarray*}


For $t \in {\cal D}_{rf}$, by Proposition \ref{bias0rx} and Remark \ref{xicon},
\begin{eqnarray*} \label {BV14a}
E \bar{g}(t; h_n)= g(t)  + \frac{\mu_{2,1} h_n^2}{2} D_2 \left ( \frac{g(t)}{  q(t)} \right ) +  \frac{\mu_{4,1} h_n^4}{24} D_4 \left ( \frac{g(t)}{ q^2(t)} \right )+ o(h^4_n).
\end{eqnarray*}
Similarly,
\begin{eqnarray}  \label {VB22ca}
E \bar{f}(t; h_n)=   f(t) + \frac{\mu_{2,1} h_n^2}{2} D_2 \left ( \frac{ f(t)}{  q(t)} \right ) +\frac{\mu_{4,1} h_n^4}{24} D_4 \left ( \frac{ f(t)}{ q^2(t)} \right )+ o(h^4_n) .
\end{eqnarray}
Hence
\begin{eqnarray}\label {32L}
&& f(t) E \bar{g}(t; h_n)-g(t)E \bar{f}(t; h_n) \nonumber\\
&=&\frac{f(t)\mu_{2,1} h_n^2}{2 }  \left \{D_2 \left ( \frac{g(t)}{  q(t)} \right ) -  r(t)  D_2 \left ( \frac{ f(t)}{  q(t)} \right ) \right \} \nonumber\\
&& +\frac{f(t)\mu_{4,1} h^4_n}{24 }\left \{ D_4 \left ( \frac{ g(t)}{  q^2(t)}\right ) -r(t) D_4 \left ( \frac{ f(t)}{q^2(t)} \right )  \right \}+o(h_n^4).
\end{eqnarray}
We next show that the term with $h^2_n$ is zero. Let $\phi(t)=f(t)/ q(t)$. Then
\begin{eqnarray*}
D_2 \left ( \frac{g(t)}{ q(t)} \right ) -  r(t)  D_2 \left ( \frac{ f(t)}{ q(t)} \right )&=& D_2 \left ( r(t) \phi(t) \right )- r(t)D_2 (\phi(t)) \nonumber\\
&=& r''(t) \phi(t) +2 r'(t)\phi'(t) \nonumber\\
&=& \frac{\left (  r'(t) \phi^2(t) \right )'}{\phi(t)} .
\end{eqnarray*}
Since $\phi(t)= f(t)/q(t)= |r'(t)|^{-1/2}$. Then $r'(t)\phi^2(t)=-1$ or $1$. Hence
\begin{eqnarray}\label {VB27ar}
D_2 \left ( \frac{g(t)}{  q(t)} \right ) -  r(t)  D_2 \left ( \frac{ f(t)}{ q(t)} \right ) = 0.
\end{eqnarray}
By (\ref{32L}) and (\ref{VB27ar}), for $\theta (t)$ defined in (\ref{VB12a}),
\begin{eqnarray}\label {37svb}
 f(t) E \bar{g}(t; h_n)-g(t)E \bar{f}(t; h_n)= f^2(t) \theta(t)h^4_n +o(h^4_n).
\end{eqnarray}
Since $E \bar{f}(t; h_n)=f(t)+o(1)$ by (\ref{VB22ca}), then by (\ref{29L}) and (\ref{37svb}),
\begin{eqnarray}\label {VB21a}
\frac{E \bar{g}(t; h_n)}{E \bar{f}(t; h_n)} = r(t) +    \theta(t) h^4_n +o(h^4_n).
\end{eqnarray}
\textbf{Step 2}. We estimate $I_1$ in (\ref{I_1}). Denote $W_i=(t-X_i) \alpha(q(X_i))/h_n$ and let $F_i=K(W_i)\alpha(q(X_i))$. By Propositons \ref{bias0rx} and \ref{bias0}, 
\begin{eqnarray}\label {30cb}
I_1 &=& \frac{1}{n^2 h_n^2} E \sum_{i \neq j} F_iY_i \left \{  F_j - EF_j \right \}  +\frac{1}{n^2 h_n^2} E \sum_{i=1}^n F_iY_i \left \{ F_i - EF_i  \right \} \nonumber \\
&=& \frac{1}{n h_n^2} E \left( F^2_1Y_1 \right)- \frac{1}{n h_n} E \left(F_1Y_1\right)\frac{1}{h_n}E F_1   \nonumber \\
&=&\frac{\sqrt{q(t)}g(t)\mu_{0,2}}{n h_n} +O\left ( \frac{1}{n}\right ).
\end{eqnarray}
\textbf{Step 3}. We estimate $I_2$ in (\ref{I_2}). By Taylor expansion of the function $\gamma(y)=(1+y)^{-1}$,
\begin{eqnarray}\label {barfrec}
\frac{1}{\bar{f}(t; h_n)}&=&\frac{1}{E\bar{f}(t; h_n)\left \{1+\bar{f}(t; h_n)/E\bar{f}(t; h_n)-1 \right\}}\nonumber\\
&=&\frac{1}{E\bar{f}(t; h_n)}\left (  1+\gamma'(\xi_t)\frac{\bar{f}(t; h_n)-E\bar{f}(t; h_n)}{E\bar{f}(t; h_n)}\right)
\end{eqnarray}
where $\xi_t$ is between $0$ and $\bar{f}(t; h_n)/E\bar{f}(t; h_n)-1$. Then
\begin{eqnarray}\label {40vbb}
I_2 &=& \frac{1}{E \bar{f}(t; h_n)}E \left \{  \bar{g}(t; h_n)  \left (\bar{f}(t; h_n) - E \bar{f}(t; h_n)\right)^2 \right \}  \nonumber \\
&& +\frac{1}{(E \bar{f}(t; h_n))^2}E \left \{ \gamma'(\xi_t) \bar{g}(t; h_n)  \left (\bar{f}(t; h_n) - E \bar{f}(t; h_n)\right)^3 \right \} \nonumber\\
&:=& \frac{1}{E \bar{f}(t; h_n)} I_{2,1}+\frac{1}{(E \bar{f}(t; h_n))^2}I_{2,2}.
\end{eqnarray}
Similar to the estimate in (\ref{30cb}),  
\begin{eqnarray} \label {37cba}
I_{2,1}&= &\frac{1 }{n^3 h^3_n } E \sum_{i\neq j}  F_i Y_i    \left \{F_j- EF_j\right \}^2   +\frac{1 }{n^3 h^3_n } E \sum_{i=1}^n F_iY_i    \left \{F_i - EF_i \right \}^2\nonumber\\
&=&\frac{1}{n h_n} E \bar{g}(t; h_n) \left \{  \frac{1}{h_n}E  F_1^2- h_n\left (\frac{1}{h_n}  E F_1\right )^2 \right \}+ O\left (\frac{1}{n^2 h^2_n}  \right)\nonumber\\
&=&\frac{\sqrt{ q(t)}g(t)f(t)\mu_{0,2}}{n h_n} +o\left ( \frac{1}{n h_n} \right).
\end{eqnarray}
By H\a"{o}lder's inequality,
\begin{eqnarray*}
|I_{2,2}| \leq \|\gamma'(\cdot)\|_{\infty}\left(E \bar{g}^2(t; h_n)\right)^{1/2}\left ( E \left (\bar{f}(t; h_n) - E \bar{f}(t; h_n)\right)^6\right)^{1/2}.
\end{eqnarray*}
Recall that $F_i = K(W_i) \alpha (q(X_i)))$. Then  
\begin{eqnarray}\label {Ebarf}
&& E \left (\bar{f}(t; h_n) - E \bar{f}(t; h_n)\right)^6\nonumber\\
&=& \frac{1}{n^6 h^6_n}E \sum_{i,j,k \ are \ different}(F_i-EF_i)^2 (F_j-EF_j)^2  (F_k-EF_k)^2 \nonumber\\
&&+  \frac{1}{n^6 h^6_n}E \sum_{i \neq j} \left \{ (F_i-EF_i)^4 (F_j-EF_j)^2+ (F_i-EF_i)^3 (F_j-EF_j)^3\right\}  \nonumber\\
&&+ \frac{1}{n^6 h^6_n}E \sum_{i=1}^n (F_i-EF_i)^6 \nonumber\\
&=& O \left (  \frac{1}{n^3 h^3_n}\right ).
\end{eqnarray}
Hence
\begin{eqnarray}\label {40cb}
I_{2,2} = o \left ( \frac{1}{n h_n} \right ).
\end{eqnarray}
Since $E \bar{f}(t; h_n)=f(t)(1+o(1))$, then by (\ref{40vbb}), (\ref{37cba}) and (\ref{40cb}), 
\begin{eqnarray}\label {42bb}
I_2 = \frac{\sqrt{q(t)}g(t)\mu_{0,2}}{n h_n} + o\left ( \frac{1}{n h_n} \right ).
\end{eqnarray}
By (\ref{21ta}), (\ref{VB21a}), (\ref{30cb}) and (\ref{42bb}), we obtain (\ref{Thm1(1)}). 

\vskip5pt
\textbf{Step 4}. 
Now we prove (\ref{Thm1(2)}). By (\ref{barfrec}),
\begin{eqnarray*}\label {46wt}
\bar{r}(t; h_n)-r(t) &=& \frac{f(t)\bar{g}(t; h_n)- g(t)\bar{f}(t; h_n)}{f(t)\bar{f}(t; h_n)}\nonumber\\
&=&\frac{f(t)\bar{g}(t; h_n)- g(t)\bar{f}(t; h_n)}{f(t)E\bar{f}(t; h_n)  }\nonumber\\
&& +  \frac{\gamma'(\xi_t) \{f(t)\bar{g}(t; h_n)- g(t)\bar{f}(t; h_n)\}\{ \bar{f}(t; h_n)-E\bar{f}(t; h_n) \}}{f(t)(E\bar{f}(t; h_n))^2 }\nonumber\\
&:=& \frac{J_1}{f(t)E\bar{f}(t; h_n) }+ \frac{J_2}{f(t)(E\bar{f}(t; h_n))^2}.
\end{eqnarray*}
Since $E \bar{f}(t; h_n)=f(t)(1+o(1))$, then
\begin{eqnarray}\label {51svb}
E \left \{   \bar{r}(t; h_n)-r(t)  \right \}^2 = \left (\frac{E J_1^2}{f^4(t) }+  \frac{2 E (J_1 J_2)}{f^5(t) } + \frac{E J_2^2}{f^6(t) }\right)(1+o(1)).
\end{eqnarray}
Recall that  $F_i=K(W_i)\alpha(q(X_i))$. Then
\begin{eqnarray}\label {45cd}
EJ_1^2&=& \frac{1}{n^2 h_n^2}E \sum_{i \neq j} \left\{f(t) F_i Y_i-g(t) F_i \right\} \left\{f(t) F_j Y_j-g(t) F_j \right\} \nonumber\\
&& + \frac{1}{n^2 h_n^2}E \sum_{i=1}^n \left\{f(t) F_i Y_i-g(t) F_i \right\}^2.
\end{eqnarray}
By (\ref{37svb}),
\begin{eqnarray}\label {46cd}
\frac{1}{h_n}E \left\{f(t) F_1 Y_1-g(t) F_1 \right\}&=& f(t) E \bar{g}(t; h_n) -g(t) E \bar{f}(t; h_n)\nonumber\\
&=& f^2(t) \theta(t)h^4_n +o(h^4_n).
\end{eqnarray}
Recall that $m(t)=E(Y_1^2|X_1=t)$ and $\sigma^2(t)=m(t)-r^2(t)$. By Proposition \ref{bias0},
\begin{eqnarray}\label {47cd}
&&\frac{1}{h_n}E \left\{f(t) F_1 Y_1-g(t) F_1 \right\}^2 \nonumber\\
&=& \frac{f^2(t)}{h_n}E  (F_1 Y_1)^2  - \frac{2 f(t)g(t)}{h_n}E F^2_1Y_1  +\frac{g^2(t)}{h_n} E F_1^2 \nonumber\\
&=& f^2(t) \sqrt{q(t)} m(t)f(t) \mu_{0,2}- 2 f(t)g(t) \sqrt{q(t)}g(t)\mu_{0,2}\nonumber\\
&&  + g^2(t) \sqrt{q(t)} f(t)\mu_{0,2}+O(h_n^2)\nonumber\\
&=&  f^3(t) \sqrt{q(t)}   \sigma^2(t) \mu_{0,2}+O(h^2_n). 
\end{eqnarray}
Applying (\ref{46cd}) and (\ref{47cd}) to (\ref{45cd}), we have
\begin{eqnarray} \label {55svb}
EJ_1^2=  \frac{f^{7/2}(t)|r'(t)|^{1/4}  \sigma^2(t) \mu_{0,2}}{n h_n} +f^4(t)\theta^2(t)h^8_n +o(h^8_n)+o\left ( \frac{1}{nh_n}\right).
\end{eqnarray}
By H\a"{o}lder's inequality,
\begin{eqnarray}\label {52wb}
E J^2_2 \leq \|\gamma'(\cdot)\|_{\infty}^2 (E J_1^4)^{1/2} \left ( E \left\{ \bar{f}(t; h_n)-E\bar{f}(t; h_n) \right\}^4\right)^{1/2}.
\end{eqnarray}
Similar to the estimate in (\ref{45cd})-(\ref{47cd}),
\begin{eqnarray}\label {53wb}
&&E J_1^4=\frac{1}{n^4 h_n^4}E \sum_{i_1, i_2, i_3, i_4 \ are \ different} \prod_{k=1}^4 \left\{f(t) F_{i_k} Y_{i_k}-g(t) F_{i_k} \right\}\nonumber\\
&& + \frac{1}{n^4 h_n^4}E\sum_{i_1, i_2, i_3 \ are \ different} \left\{f(t) F_{i_1} Y_{i_1}-g(t) F_{i_1} \right\}^2\prod_{k=2}^3 \left\{f(t) F_{i_k} Y_{i_k}-g(t) F_{i_k} \right\} \nonumber\\
&&+\frac{1}{n^4 h_n^4}E \sum_{i \neq j}\left\{f(t) F_i Y_i-g(t) F_i \right\}^2 \left\{f(t) F_j Y_j-g(t) F_j \right\}^2\nonumber\\
&&+ \frac{1}{n^4 h_n^4}E \sum_{i \neq j} \left\{f(t) F_i Y_i-g(t) F_i \right\}^3 \left\{f(t) F_j Y_j-g(t) F_j \right\}\nonumber\\
&&+ \frac{1}{n^4 h_n^4}E \sum_{i =1}^n \left\{f(t) F_i Y_i-g(t) F_i \right\}^4 \nonumber\\
&=&f^8(t) \theta^4 (t)h^{16}_n +o(h^{16}_n) + O \left (\frac{h_n^8}{n h_n}  \right )+O \left ( \frac{1}{n^2 h_n^2} \right ).
\end{eqnarray}
Similar to (\ref{Ebarf}),
\begin{eqnarray}\label {54wb}
  E \left (\bar{f}(t; h_n) - E \bar{f}(t; h_n)\right)^4= O \left (  \frac{1}{n^2 h^2_n}\right ).
\end{eqnarray}
Applying (\ref{53wb}) and (\ref{54wb}) to (\ref{52wb}), we have
\begin{eqnarray} \label {55wa}
E J_2^2 = o(h^8_n)+ o \left (\frac{1}{n h_n}  \right ).
\end{eqnarray}
By (\ref{55svb}) and (\ref{55wa}), and by H\a"{o}lder's inequality, $|E(J_1J_2)| = o(h^8_n)+o(1/(n h_n))$. Then by (\ref{51svb}), (\ref{55svb}) and (\ref{55wa}), we obtain (\ref{Thm1(2)}). 
\end{proof}
\bigskip


\begin{proof}[Proof of Theorem \ref{central1}]
 By (\ref{Thm1(1)}) and (\ref{Thm1(2)}),
\begin{eqnarray*}
E \bar{r}(t; h_n) - r(t)=  \theta(t) h_n^4   +o(h_n^4)+o\left (\frac{1}{n h_n} \right ) 
\end{eqnarray*}
and
\begin{eqnarray*}
Var\left ( \bar{r}(t; h_n)- r(t) \right )= \frac{\mu_{0,2}|r'(t)|^{1/4}\sigma^2(t)}{ nh_n \sqrt{f(t)}}  +o(h_n^8)+o\left (\frac{1}{n h_n}  \right).
\end{eqnarray*}
If $h_n^4 \sqrt{n h_n} \rightarrow \lambda$, then by the Lindeberg's central limit theorem, 
\begin{eqnarray*}
\sqrt{n h_n} \left \{\bar{r}(t; h_n)-r(t) \right \}\xrightarrow {D}N \left (\lambda \theta(t),  \frac{\mu_{0,2}|r'(t)|^{1/4}\sigma^2(t)}{\sqrt{f(t)}} \right ).
\end{eqnarray*}
\end{proof}


To prove Theorem \ref{Thm2}, we first establish the following two lemmas. Note that $\hat{q}(t; h_n)$ is a function of all observations $X_1,..., X_n$.


\begin{lemma}\label {Lemma1}
Under the condition of Theorem \ref{Thm2},  for any integer $a \geq 1$, 
\begin{eqnarray*}\label {46vbaa}
&&  E \left \{\left (\hat{q}^4 (X_1; h_{1,n})- q^4 (X_1 ) \right) \big |X_1,...,X_a \right \} \nonumber\\
&=&  \Psi(X_1, h_{1,n}) h_{1,n}^2+   \frac{\phi(X_1)}{n h^3_{1,n}} +O\left (  \frac{1}{n h^2_{1,n}}\right)+o(h^6_{2,n})
\end{eqnarray*}
for some functions $\Psi$ and $\phi$, where $\Psi(x, h_{1,n})=b(x)+d(x)h^2_{1,n}$ for some functions $b(x)$ and $d(x)$.
\end{lemma}
\begin{proof}
Fix $X_1=x_1$. Since $\hat{r}(x_1; h_{1,n})=\hat{g}(x_1;h_{1,n})/\hat{f}(x_1; h_{1,n}) $, then by (\ref{hatq}),
\begin{eqnarray*}
\hat{q}^4(x_1; h_{1,n})&=& \left \{\hat{f}^2(x_1; h_{1,n})\hat{r}'(x_1; h_{1,n}) \right \}^2 \nonumber\\
&=& \left \{  \hat{f}(x_1; h_{1,n})\hat{g}'(x_1; h_{1,n}) - \hat{g}(x_1; h_{1,n}) \hat{f}'(x_1; h_{1,n})\right \}^2,
\end{eqnarray*}
where  
\begin{eqnarray*} 
 \hat{g}'(x_1; h_{1,n}) :=\frac{d\hat{g}(x; h_{1,n}) }{dx}\Big|_{x=x_1} = \frac{1}{n h^2_{1,n}}\sum_{i=1}^n K'\left ( \frac{x_1-X_i}{h_{1,n}}    \right )     Y_i, \label {VB58aw} \\
 \hat{f}'(t; h_{1,n}):=\frac{d\hat{f}(x_1; h_{1,n}) }{dt}\Big|_{x=x_1} = \frac{1}{n h^2_{1,n}}\sum_{i=1}^n K'\left ( \frac{x_1-X_i}{h_{1,n}}    \right ) . \label {VB59aw}
\end{eqnarray*}
For $m \in [1,n]$, denote
\begin{eqnarray}\label {Nts}
U^{(m)}_i &=& (X_m-X_i)/h_{1,n} \nonumber\\
G^{(m)}_{i,j} &=&K(U^{(m)}_i) K'(U^{(m)}_j)Y_j  - K(U^{(m)}_i)Y_i K'(U^{(m)}_j) \nonumber\\
H^{(m)}_{i,j,k,l}&=&G^{(m)}_{i,j}G^{(m)}_{k,l} - h^6_{1,n} q^4(X_m)
\end{eqnarray}
Then $G^{(m)}_{i,i}=0$ for any $m,i \in [1,n]$, and
\begin{eqnarray}\label {64st}
\hat{q}^4(x_1; h_{1,n})-q^4(x_1) &=&\frac{1}{n^4 h^6_{1,n}}\sum_{ i,j,k, l }  \left  (G^{(1)}_{i,j}G^{(1)}_{k,l}  -h^6_{1,n}q^4(x_1)\right)\nonumber\\
&=& \frac{1}{n^4 h^6_{1,n}}\sum_{ i,j,k, l }  H^{(1)}_{i,j,k,l} .
\end{eqnarray}
Denote $E^*$ as the expectation of $\{ X_{a+1},..., X_n\}$. Write
\begin{eqnarray}\label {E*}
 E^*[ \hat{q}^4(x_1; h_{1,n})-q^4(x_1)] =J_1+J_2+J_3,
\end{eqnarray}
where
\begin{eqnarray}
J_1&=&E^* \Bigg(\frac{1}{n^4 h^6_{1,n} } \sum_{1,  i, j, k, l \;\;are \;\;different  }  H^{(1)}_{i,j, k, l}  \Bigg),  \nonumber \\
J_2&=& E^* \Bigg(\frac{1}{n^4 h^6_{1,n} } \sum_{  exactly \ two \ of \ 1,   i, j, k, l \ are \ equal   }  H^{(1)}_{i,j, k , l} \Bigg),  \nonumber\\
J_3&=&  E^* \Bigg(\frac{1}{n^4 h^6_{1,n} } \sum_{\substack{  other \  cases  }}  H^{(1)}_{i,j, k , l} \Bigg). \nonumber
\end{eqnarray}
First we estimate $E^* G^{(1)}_{i,j}/h^3_{1,n}$ for $i \neq  j$ and $i,j \in [a+1, n]$. Denote $w_j=(x_1-x_j)/h_{1,n}$. Then
\begin{eqnarray}\label {63vbag}
\frac{1}{h_{1,n}^3} E^* G^{(1)}_{i,j}=  \frac{1}{h_{1,n}^3} \int K'(w_3) \int K(w_2)  (r(x_3)-r(x_2))f(x_2)dx_2 f(x_3)dx_3.
\end{eqnarray}
Since $g(x_2)=r(x_2)f(x_2)$, by Proposition \ref{bias0rx} with $\xi(s)=1$,
\begin{eqnarray}\label {64vbag}
&& \frac{1}{h_{1,n}}\int K(w_2)  (r(x_3)-r(x_2))f(x_2)dx_2 \nonumber\\
&&=r(x_3) \left \{f(x_1) +\frac{f''(x_1)\mu_{2,1}}{2} h^2_{1,n} \right \}\nonumber\\
&& -\left \{ g(x_1)+\frac{g''(x_1)\mu_{2,1}}{2}h^2_{1,n} \right\}+O(h_{1,n}^4). \qquad
\end{eqnarray}
By (\ref{63vbag}) and (\ref{64vbag}),
\begin{eqnarray}\label {65vbag}
 &&\frac{1}{h_{1,n}^3} E^* G^{(1)}_{i,j}=\left \{f(x_1) +\frac{f''(x_1)\mu_{2,1}}{2} h^2_{1,n} \right \}\frac{1}{h_{1,n}^2}\int K'(w_3)g(x_3)dx_3 \nonumber\\
&& - \left \{ g(x_1)+\frac{g''(x_1)\mu_{2,1}}{2}h^2_{1,n} +O(h_{1,n}^4)\right\}\frac{1}{h_{1,n}^2}\int K'(w_3)f(x_3)dx_3 . \qquad
\end{eqnarray}
Note that $\int K'(w)w^{2k}dw=0$ for $k =0, 1,\cdots$, $\int K'(w)wdw=-1$ and $\int K'(w)w^3 dw=-3\mu_{2,1}$. Then by Taylor expansion,
\begin{eqnarray}\label {66vbag}
&& \frac{1}{h_{1,n}^2}\int K'(w_3)g(x_3)dx_3 \nonumber\\
&=& \frac{1}{h_{1,n}}\int K'(w_3) \bigg \{g(x_1)-g'(x_1)w_3 h_{1,n} +\frac{g''(x_1)}{2}w_3^2 h_{1,n}^2\nonumber\\
&&-\frac{g'''(x_1)}{6} w^3_{3}h_{1,n}^3 +\frac{g^{(4)}(x_1)}{24}w_3^4 h_{1,n}^4 \bigg \}dw_3 +O(h_{1,n}^4)   \nonumber\\
&=&g'(x_1)+\frac{g'''(x_1)\mu_{2,1}}{2}h_{1,n}^2 +O(h_{1,n}^4).
\end{eqnarray}
Similarly,
\begin{eqnarray}\label{67vbag}
\frac{1}{h_{1,n}^2}\int K'(w_3)f(x_3)dx_3= f'(x_1)+\frac{f'''(x_1)\mu_{2,1}}{2}h_{1,n}^2 +O(h_{1,n}^4). 
\end{eqnarray}
Applying (\ref{66vbag}) and (\ref{67vbag}) to (\ref{65vbag}), we have, for some function $b_0$,
\begin{eqnarray*}\label {64vbr}
\frac{1}{h_{1,n}^3} E^* G^{(1)}_{i,j}&=& f(x_1)g'(x_1)-g(x_1)f'(x_1)+b_0(x_1)h^2_{1,n}+O(h_{1,n}^4). \qquad
\end{eqnarray*}
Since $q^4(x_1)=( f(x_1)g'(x_1)-g(x_1)f'(x_1))^2$, then for some function $b$,
\begin{eqnarray}\label {64tt}
\left ( \frac{1}{ h^3_{1,n}}E^* G^{(1)}_{i,j}  \right)^2 - q^4(x_1) =  b (x_1)h^2_{1,n}+ O(h^4_{1,n}). 
\end{eqnarray}
If we apply the above Taylor expansion further, $O(h_{1,n}^4)$ in (\ref{64tt}) can be expressed as $d(x_1)h^4_{1,n}+O(h^6_{1,n})$ for some funtion $d$. Then  \\
\begin{eqnarray} \label {65vbr}
 \left ( \frac{1}{ h^3_{1,n}}E^*G^{(1)}_{i,j}  \right)^2 - q^4(x_1)=   \Psi(x_1, h_{1,n})h^2_{1,n} +O(h^6_{1,n}) 
\end{eqnarray}
where $\Psi(x_1, h_{1,n})=b(x_1) +d(x_1)h^2_{1,n}$. 

Next we estimate $J_1$ in (\ref{E*}) with three cases. Case 1: $i,j,k,l \in [a+1, n]$. By (\ref{65vbr}),
\begin{eqnarray*} \label {57vt}
&&E^* \Bigg(\frac{1}{n^4 h^6_{1,n} } \sum_{\substack{ i,j,k, l \in [a+1, n]  \\ 1, i, j, k, l \;\;are \;\;different  }}  H^{(1)}_{i,j, k, l}  \Bigg) \nonumber\\  
&&= \Psi(x_1, h_{1,n})   h^2_{1,n}+O(h^6_{1,n})    + O \left (  \frac{1}{n}\right).\qquad
\end{eqnarray*}
Case 2: One of $i,j,k,l$ is in $[2,a]$. If $j \in [2,a]$, for example, then
\begin{eqnarray}\label {72vt}
&&  \frac{1}{n^4 h^6_{1,n}} E^* \Bigg(\sum_{\substack{i,k,l \in [a+1, n], \ j  \in [2,a]\\ 1, i,k, l \ are \ different}}H^{(1)}_{i,j,k,l} \Bigg) \nonumber\\
&=&  \frac{   K'(\frac{x_1-X_j}{h_{1,n}})r(X_j)E^* K(U^{(1)}_i) -K'(\frac{x_1-X_j}{h_{1,n}})E^* \left \{K(U^{(1)}_i)r(X_i) \right \}}{n h^3_{1,n}}\times \frac{E^* G^{(1)}_{k,l}}{h^3_{1,n}} -\frac{q^4(x_1)}{n}\nonumber
\end{eqnarray}
where $E^* K(U^{(1)}_i)/h_{1,n}=E^*\hat{f}(x_1; h_{1,n})$ and $E^* \left \{K(U^{(1)}_i)r(X_i) \right\}/h_{1,n}=E^*\hat{g}(x_1, h_{1,n})$. Hence
\begin{eqnarray*}
\frac{1}{n^4 h^6_{1,n}} E^* \Bigg(\sum_{\substack{i,k,l \in [a+1, n], \ j \in [2,a]\\ i,k, l \ are \ different}}H^{(1)}_{1, i,j,k,l} \Bigg)= O\left (  \frac{1}{n h^2_{1,n}}\right) .
\end{eqnarray*}
If $i,k$ or $l$ is in $[2,a]$, the results are similar. Hence
\begin{eqnarray*} \label {65ww}
E^* \Bigg(\frac{1}{n^4 h^6_{1,n} } \sum_{\substack{one \ of \  i,j,k, l \ in \ [2, a]  \\1,  i, j, k, l \;\;are \;\;different  }}  H^{(1)}_{i,j, k, l}  \Bigg)  =O\left (  \frac{1}{n h^2_{1,n}}\right).
\end{eqnarray*}
Case 3: Two or more of $i,j,k,l$ are in $[2,a]$. The result is $O(1/(n^2 h^4_{1,n}))=o(h_{2,n}^{10})$. Combining the three cases, we have
\begin{eqnarray}\label {E*1}
J_1 =  \Psi(x_1, h_{1,n}) h_{1,n}^2+o(h^6_{2,n}) +O\left (  \frac{1}{n h^2_{1,n}}\right). 
\end{eqnarray}
Now we estimate $J_2$ with two cases. Case 1: Two of $i,j,k,l$ are equal, and $i,j,k,l \in [a+1, n]$. For example, if $j=l$ and $i,j,k \in [a+1, n]$ are different, then
\begin{eqnarray*}\label {69st}
&& \frac{1}{n^4 h^6_{1,n}}E^* \Bigg(\sum_{\substack{i,j,k \in [a+1, n]\\ i , j, k \ are \ different}}H^{(1)}_{i,j,k,j}\Bigg)\nonumber\\
&=& \frac{1}{n h^6_{1,n}} E^* \left \{K(U^{(1)}_i)K'(U^{(1)}_j) (Y_j-Y_i)K(U^{(1)}_k)K'(U^{(1)}_j) (Y_j-Y_k)\right\}\nonumber\\
&&- \frac{q(x_1)}{n}+O\left( \frac{1}{n}\right)\qquad
\end{eqnarray*}
where the first term on the right-hand side turns out to be a function of $x_1$ multiplied by $1/(n h^3_{1,n})$. The other situations have similar results. Case 2: Two of $i,j,k,l$ are equal, and some $i,j,k,l$ are in $[1,a]$. The result is $O(1/(n^2 h^6_{1,n}))=o(h_{2,n}^{8})$. Then
\begin{eqnarray}\label {E*2}
J_2 =   \frac{\phi(x_1)}{n h^3_{1,n}}+ o(h_{2,n}^{8}).
\end{eqnarray}
for some function $\phi$. Also
\begin{eqnarray}\label {E*3}
J_3=O\left ( \frac{1}{n^2 h^4_{1,n}} \right)=o(h_{2,n}^{10}).
\end{eqnarray}
Applying (\ref{E*1}), (\ref{E*2}) and (\ref{E*3}) to (\ref{E*}), we obtain the lemma. 
\end{proof}
\bigskip


\begin{lemma}\label {Lemma2}
Denote $W_i=(t-X_i)\alpha(q(X_i))/h_{2,n}$ and let $L(w)= K(w)+w K'(w)$. Assume that $f'/f$ and $r''/r'$ are bounded in a neighborhood of $t$. 
 For any function $H$ with bounded and continuous second order derivtive,   
\begin{eqnarray*}
\frac{1}{h_{2,n}} E \left \{L(W_1)H(X_1)   \right \} = O(h^{2}_{2,n}).
\end{eqnarray*}
\end{lemma}
\begin{proof}
We have
\begin{eqnarray*}\label {98Vt}
 \frac{1}{h_{2,n}} E \left \{L(W_1)H(X_1)  \right \} =\frac{1}{h_{2,n}} \int L \left ( (t-s)\alpha(q(s))/h_{2,n} \right ) H(s) ds. 
\end{eqnarray*}
Let $U_t(v)=v\alpha(q(t-v))$. By a similar argument as in Remark \ref{xicon},  $U'_t(v)=\alpha(q(t-v))-v\alpha'(q(t-v))q'(t-v)>0$ in a small neighborhood of $0$. Then $U_t(v)$ is invertible in a small neighborhood of $0$.  Denote $V_t(u)$ as the inverse function of $U_t(v)$. Let $zh_{2,n}= (t-s)\alpha (q(s))= (t-s)\alpha (q(t-(t-s)))$. Then $z h_{2,n}=U_t(t-s)$ and hence $t-s = V_t(zh_{2,n})$.  The change of variables from $s$ to $z$ gives
\begin{eqnarray*}\label{G}
\frac{1}{h_{2,n}} E L(W_1) H(X_1)&=& - \int L \left ( z \right ) H \left (t- V_t(zh_{2,n})\right) \frac{dV_t}{du}_{u=zh_{2,n}}dz \nonumber\\
&:=& - \int L \left ( z \right ) G_t(zh_{2,n})dz. 
\end{eqnarray*}
By Taylor expansion, $G_t(zh_{2,n})=G_t(0)+G_t'(0)zh_{2,n}+ G_t''(\tau)(z^2h_{2,n}^2)/2$ where $\tau$ is between $0$ and $zh_{2,n}$. By the condition on $H$,  $G_t$  has bounded and continuous second order derivtive. Also notice that $\int L(z)dz=0$, $\int z L(z)dz=0$. Then
\begin{eqnarray*}\label{Lq'}
\frac{1}{h_{2,n}} E\left \{ L(W_1) H (X_1)\right\} = O(h^2_{2,n}).
\end{eqnarray*}
\end{proof}


\begin{proof}[Proof of Theorem \ref{Thm2}]

By (\ref{VB16ag}) with $p=1$ and $z=\hat{f}(t; h_{1,n}, h_{2,n})/E \hat{f}(t; h_{1,n}, h_{2,n})$,
\begin{eqnarray}  \label {75vbt} 
E \hat{r} (t; h_{1,n}, h_{2,n}) = \frac{E \hat{g}(t; h_{1,n}, h_{2,n})}{E \hat{f}(t; h_{1,n}, h_{2,n})}   +\frac{-M_1+M_2}{(E \hat{f} (t; h_{1,n}, h_{2,n}))^2}  
\end{eqnarray}
where
\begin{eqnarray}\label {I1}
M_1=E\left \{\hat{g}(t; h_{1,n}, h_{2,n}) \left (\hat{f}(t; h_{1,n}, h_{2,n})-E \hat{f}(t; h_{1,n}, h_{2,n})\right) \right\}
\end{eqnarray}
and
\begin{eqnarray}\label {I2}
M_2=E \left \{\hat{r}(t; h_{1,n}, h_{2,n}) \left (\hat{f}(t; h_{1,n}, h_{2,n})-E \hat{f}(t; h_{1,n}, h_{2,n})\right)^2 \right\}.
\end{eqnarray}
\textbf{Step 1}. We estimate $E \hat{g}(t; h_{1,n}, h_{2,n})/E\hat{f}(t; h_{1,n}, h_{2,n})$. Define the function $\beta(y)=\alpha( y^{1/4})$ for $y> 0$. By Taylor expansion, 
\begin{eqnarray}\label {65vbc}
 \alpha (\hat{q}(X_i; h_{1,n})) -\alpha(q(X_i))&= & \beta(\hat{q}^4(X_i; h_{1,n}))-\beta(q^4(X_i))\nonumber\\
 &=& \beta'(q^4(X_i)) \left \{ \hat{q}^4(X_i; h_{1,n}) - q^4(X_i)   \right \}  \nonumber\\
&&+\frac{\beta''(\hat{\nu_i})}{2} \left \{\hat{q}^4(X_i; h_{1,n}) -q^4(X_i) \right \}^2 \qquad
\end{eqnarray}
where $\hat{\nu_i}$ is between $\hat{q}^4(X_i; h_{1,n})$ and $q^4(X_i)$. Denote $W_i=(t-X_i) \alpha (q(X_i))/h_{2,n}$.  By Taylor expansion and (\ref{65vbc}),
\begin{eqnarray} \label {66vbc}
&& K \left ( \frac{t-X_i}{h_{2,n}}\alpha (\hat{q}(X_i; h_{1,n})) \right ) \nonumber\\
&=&K(W_i)+\sum_{k=1}^2 \frac{K^{(k)} \left (W_i\right) }{k!}   \left (\frac{t-X_i}{h_{2,n}}\left \{ \beta(\hat{q}^4(X_i; h_{1,n}))-\beta(q^4(X_i)) \right \}\right )^k\nonumber\\
&& + \frac{K'''   (\hat{\xi}_i ) }{6}   \left (\frac{t-X_i}{h_{2,n}}\left \{ \beta(\hat{q}^4(X_i; h_{1,n}))-\beta(q^4(X_i)) \right \}\right )^3 \nonumber\\
&=& K  \left (W_i\right)+K'  \left (W_i\right) \frac{t-X_i}{h_{2,n}}\beta'(q^4(X_i))\left \{ \hat{q}^4(X_i; h_{1,n}) - q^4(X_i)   \right \} + \hat{\delta}_{i}
\end{eqnarray}
where $\hat{\xi}_i$ is between $ (t-X_i)\alpha( \hat{q}(X_i; h_{1,n}) )/h_{2,n}$ and $W_i$, and
\begin{eqnarray}\label{delta}
\hat{\delta}_{i}&=&K'  \left (W_i\right) \frac{t-X_i}{h_{2,n}}\Bigg \{  \frac{1}{2}\beta''(\hat{\nu}_i) \left \{ \hat{q}^4(X_i; h_{1,n}) - q^4(X_i)   \right \}^2    \Bigg \} \nonumber\\
&&+  \frac{K''\left (W_i\right) }{2}   \left (\frac{t-X_i}{h_{2,n}}\left \{  \beta(\hat{q}^4(X_i; h_{1,n}))-\beta(q^4(X_i))\right \}\right )^2\nonumber\\
&& + \frac{K'''   (\hat{\xi}_i ) }{6}   \left (\frac{t-X_i}{h_{2,n}}\left \{ \beta(\hat{q}^4(X_i; h_{1,n}))-\beta(q^4(X_i))\right \}\right )^3.
\end{eqnarray}
Let $L(w)= K(w)+w K'(w)$. Then  by (\ref{65vbc}) and (\ref{66vbc}),
\begin{eqnarray}\label {71vbd}
&& K \left ( \frac{t-X_i}{h_{2,n}}\alpha (\hat{q}(X_i; h_{1,n})) \right )\alpha (\hat{q}(X_i; h_{1,n})) \nonumber\\
&=&K  \left (W_i\right) \alpha( q(X_i))+ L(W_i) \beta'(q^4(X_i)) \left \{\hat{q}^4(X_i; h_{1,n}) -q^4(X_i) \right \}+ \hat{\eta}_{i} \qquad
\end{eqnarray}
where
\begin{eqnarray}\label {83vbs}
\hat{\eta}_{i}&=& K(W_i)  \left \{  \frac{1}{2}\beta''(\hat{\nu}_i) \left \{ \hat{q}^4(X_i; h_{1,n}) - q^4(X_i)   \right \}^2   \right \}  \nonumber\\
&& + K'  \left (W_i\right) \frac{t-X_i}{h_{2,n}}\beta'(q^4(X_i))\left \{ \hat{q}^4(X_i; h_{1,n}) - q^4(X_i)   \right \} \nonumber\\
&&\times\left \{\beta(\hat{q}^4(X_i; h_{1,n}))-\beta(q^4(X_i))\right \} \nonumber\\
&&+\hat{\delta}_i  \alpha (\hat{q}(X_i; h_{1,n})).
\end{eqnarray}
Similar to the proof of Lemma \ref{Lemma1},  by (\ref{delta}), (\ref{83vbs}),  and using the second equation of (\ref{65vbc}),
\begin{eqnarray}\label {81ce}
\frac{1}{ h_{2,n}} E \left\{ E  (\hat{\eta}_{1}Y_1  |X_1)\right\}= o(h^4_{2,n}) . 
\end{eqnarray}
Let
\begin{eqnarray} \label {lgxi}
\hat{\theta}_i =  L(W_i)  \beta'(q^4(X_i))\left \{\hat{q}^4(X_i; h_{1,n}) -q^4(X_i) \right \}+ \hat{\eta}_{i}  .
\end{eqnarray}
By (\ref{71vbd}),
\begin{eqnarray}\label {94Eg}
  \hat{g}(t; h_{1,n}, h_{2,n})= \bar{g}(t; h_{2,n})+ \frac{1}{n h_{2,n}} \sum_{i=1}^n \hat{\theta}_i Y_i 
\end{eqnarray}
and
\begin{eqnarray}\label {123vbaw}
  \hat{f}(t; h_{1,n}, h_{2,n})=  \bar{f}(t; h_{2,n})+ \frac{1}{n h_{2,n}} \sum_{i=1}^n \hat{\theta}_i .
\end{eqnarray}
By Lemma \ref{Lemma1} and Lemma \ref{Lemma2},\\
\begin{eqnarray} \label {91pa}
&& \frac{1}{h_{2,n}}E \left \{  L(W_1) \beta'(q^4(X_1))Y_1 E\left( \left \{\hat{q}^4(X_1; h_{1,n}) -q^4(X_1) \right \}|X_1 \right) \right\}\nonumber\\
&=& \frac{1}{h_{2,n}}E \left \{L(W_1)\beta'(q^4(X_1))  r(X_1) \left (\Psi(X_1, h_{1,n})  h^2_{1,n} +\frac{\phi(X_1)}{n h^3_{1,n}}+ O\left ( \frac{1}{n h^2_{1,n}} \right)+o(h^6_{2,n})\right) \right\}  \nonumber\\
&=& O (h^2_{1,n}h^2_{2,n})+\frac{O(h^2_{2,n})}{n h^3_{1,n}}\nonumber\\
&&+\left \{ O \left (  \frac{1}{n h^2_{1,n}}\right )+o(h^6_{1,n})\right\}\frac{1}{h_{2,n}}E \left |  L(W_1)\beta'(q^4(X_1))  r(X_1)\right|\nonumber\\
&=& o(h^4_{2,n}).
\end{eqnarray}
Applying (\ref{81ce}) and (\ref{91pa}) to (\ref{lgxi}), we have
\begin{eqnarray}\label {92pa}
\frac{1}{ h_{2,n}}E   \{E   ( \hat{\theta}_1 Y_1  |X_1) \} = o(h^4_{2,n}).
\end{eqnarray}
By (\ref{94Eg}) and (\ref{92pa}),
\begin{eqnarray}\label {94vbt}
E  \hat{g}(t; h_{1,n}, h_{2,n})= E \bar{g}(t; h_{2,n})+  o(h^4_{2,n}).
\end{eqnarray}
Similarly,
\begin{eqnarray}\label {95vbt}
E \hat{f}(t; h_{1,n}, h_{2,n})= E \bar{f}(t; h_{2,n})+  o(h^4_{2,n}) .
\end{eqnarray}
By (\ref{VB21a}) in the proof of Theorem \ref{Thm01},
\begin{eqnarray}\label {96vbt}
\frac{E \hat{g}(t; h_{1,n}, h_{2,n})}{E \hat{f}(t; h_{1,n}, h_{2,n})} = \frac{E\bar{g}(t; h_{2,n})}{E\bar{f}(t; h_{2,n})}+o(h^4_{2,n})= r(t) +    \theta(t) h^4_{2,n }+   o(h^4_{2,n}).
\end{eqnarray}
\textbf{Step 2}. We estimte $M_1$ in (\ref{I1}). Applying (\ref{94vbt})  and (\ref{95vbt}) to (\ref{I1}), we have
\begin{eqnarray}\label {I1a}
M_1 = E \left\{\hat{g}(t; h_{1,n}, h_{2,n}) \hat{f}(t; h_{1,n}, h_{2,n})\right\}- E \bar{g}(t; h_{2,n}) E \bar{f}(t; h_{2,n})+o(h^4_{2,n}) .
\end{eqnarray}
By (\ref{94Eg}) and (\ref{123vbaw}),
\begin{eqnarray}\label {124vbaw}
 && E\left\{\hat{g}(t; h_{1,n}, h_{2,n}) \hat{f}(t; h_{1,n}, h_{2,n})\right\}\nonumber\\
&=& E\left\{ \bar{g}(t; h_{2,n})\bar{f}(t; h_{2,n})\right\}+ E\left\{\bar{f}(t; h_{2,n})\frac{1}{n h_{2,n}} \sum_{i=1}^n \hat{\theta}_i  Y_i\right\}\nonumber\\
&&+E \left \{\bar{g}(t; h_{2,n})\frac{1}{n h_{2,n}}  \sum_{i=1}^n \hat{\theta}_i \right\}  + E \left \{\frac{1}{n^2 h^2_{2,n}}\sum_{i=1}^n  \hat{\theta}_iY_i  \sum_{i=1}^n  \hat{\theta}_i \right\}.
\end{eqnarray}
Similar to (\ref{92pa}),
\begin{eqnarray}\label {124Ld}
&& E\left\{\bar{f}(t; h_{2,n})\frac{1}{n h_{2,n}}  \sum_{i=1}^n \hat{\theta}_i  Y_i\right\} \nonumber\\
&=&\frac{1}{n^2 h^2_{2,n}} E \sum_{i\neq j}K(W_j)\alpha(q(X_j))  \hat{\theta}_iY_i   + \frac{1}{n^2 h^2_{2,n}}E\sum_{i=1}^nK(W_i)\alpha(q(X_i))  \hat{\theta}_iY_i \nonumber\\
&=& o(h^4_{2,n}) .
\end{eqnarray}
Similarly, the third term on the right-hand side of (\ref{124vbaw}) is $o(h^4_{2,n})$. Note that
\begin{eqnarray}\label {90vt}
 E \left \{\frac{1}{n^2 h^2_{2,n}}\sum_{i=1}^n  \hat{\theta}_iY_i  \sum_{i=1}^n  \hat{\theta}_i \right\}= \frac{1}{n^2 h^2_{2,n}}E\sum_{i\neq j}\hat{\theta}_iY_i \hat{\theta}_j +  \frac{1}{n^2 h^2_{2,n}}E\sum_{i}\hat{\theta}^2_iY_i .
\end{eqnarray}
For later use, we now show that it is $o(h^8_{2,n})+o(1/(n h_{2,n}))$. By (\ref{lgxi}) and letting $\hat{\lambda}_i=L(W_i)  \beta'(q^4(X_i))\left \{\hat{q}^4(X_i; h_{1,n}) -q^4(X_i) \right \}$, we have
\begin{eqnarray}\label {101tt}
 \frac{1}{n^2 h^2_{2,n}}E\sum_{i\neq j}\hat{\theta}_iY_i \hat{\theta}_j =  \frac{1}{n^2 h^2_{2,n}}E\sum_{i\neq j}\left ( \hat{\lambda}_i Y_i + \hat{\eta}_{i}Y_i\right)  \left ( \hat{\lambda}_j   + \hat{\eta}_{j}\right).  
\end{eqnarray}
By (\ref{Nts}) and (\ref{64st}),
\begin{eqnarray}\label {102tt}
 \frac{1}{n^2 h^2_{2,n}}E\sum_{i\neq j} \hat{\lambda}_i Y_i \hat{\lambda}_j &=&  \frac{1}{n^2 h^2_{2,n}}E\sum_{i\neq j}\bigg\{ L(W_i)  \beta'(q^4(X_i))L(W_j)  \beta'(q^4(X_j)) \nonumber\\
&& \times \frac{1}{n^4 h^6_{1,n}}\sum_{k, l , p, q}  H^{(i)}_{k,l, p, q}  \frac{1}{n^4 h^6_{1,n}}\sum_{k', l' , p', q'}  H^{(j)}_{k',l', p', q'} \bigg \}.\qquad
\end{eqnarray}
Write
\begin{eqnarray}
\frac{1}{n^4 h^6_{1,n}}\sum_{k, l , p, q}  H^{(i)}_{k,l, p, q} \frac{1}{n^4 h^6_{1,n}}\sum_{k', l' , p', q'}  H^{(j)}_{k',l', p', q'}=T_1 +T_2 +T_3
\end{eqnarray}
where
\begin{eqnarray}
T_1 &=& \frac{1}{n^8 h^{12}_{1,n}}\sum_{i, j , k,l,p,q,  k', l' , p', q' \ are \ different}H^{(i)}_{k,l, p, q}   H^{(j)}_{k',l', p', q'}, \nonumber\\
T_2 &=& \frac{1}{n^8 h^{12}_{1,n}}\sum_{exactly \ two \ of \ i, j , k,l,p,q,  k', l' , p', q' \ are \ equal} H^{(i)}_{k,l, p, q}  H^{(j)}_{k',l', p', q'}, \nonumber\\
T_3 &=& \frac{1}{n^8 h^{12}_{1,n}}\sum_{Other \ cases} H^{(i)}_{k,l, p, q}  H^{(j)}_{k',l', p', q'}.
\end{eqnarray}
Similar to the estimate of $J_1$ in Lemma \ref{Lemma1}, and by Lemma \ref{Lemma2},
\begin{eqnarray}\label {105tt}
&& \frac{1}{n^2 h^2_{2,n}}E\sum_{i\neq j}\left\{ L(W_i)  \beta'(q^4(X_i))L(W_j)  \beta'(q^4(X_j)) T_1 \right\} \nonumber\\
&=&  \frac{1}{n^2 h^2_{2,n}}\sum_{i\neq j}E\left\{ L(W_i)  \beta'(q^4(X_i)) \left ( \Psi(X_i, h_{1,n}) h_{1,n}^2\right ) \right\}\nonumber\\
&& \times E\left\{ L(W_j)  \beta'(q^4(X_j)) \left ( \Psi(X_j, h_{1,n}) h_{1,n}^2  \right ) \right\}+o(h^8_{2,n})+O\left (\frac{1}{n}  \right )\nonumber\\
&=& o(h^8_{2,n})+o\left ( \frac{1}{n h_{2,n}}\right ). 
\end{eqnarray}
We consider two cases  when replacing $T_1$ by $T_2$ in the above analysis. \\
Case 1: The two equal indices are both in $\{i, k, l, p, q  \}$ or both in $\{ j, k',l', p', q'\}$. Suppose they are in $\{i, k, l, p, q  \}$. Similar to the estimates of $J_2$ and $J_1$ in Lemma \ref{Lemma1}, the result is
\begin{eqnarray} \label {106tt}
&&  \frac{1}{n^2 h^2_{2,n}}\sum_{i\neq j}E\left\{ L(W_i)  \beta'(q^4(X_i)) \left ( \frac{\phi(X_i)}{n h^3_{1,n}} \right ) \right\}\nonumber\\
&& \times E\left\{ L(W_j)  \beta'(q^4(X_j)) \left ( \Psi(X_j, h_{1,n}) h_{1,n}^2 \right ) \right\}+o(h^8_{2,n})+o\left ( \frac{1}{n h_{2,n}} \right )\nonumber\\
&=&o(h^8_{2,n})+o\left ( \frac{1}{n h_{2,n}} \right ). 
\end{eqnarray}
Case 2: One of the equal indices is in $\{i, k, l, p, q  \}$ and the other in $\{ j, k',l', p', q'\}$. Suppose that $k=k'$. In this case we first fix $X_k$ and take conditional expectation of the other variables, and then take expectation of $X_k$. The result is $O(h^4_{2,n}/(n h^3_{1,n}))=o(h^8_{2,n})$. Combining the two cases, we conclude
\begin{eqnarray*}
\frac{1}{n^2 h^2_{2,n}}E\sum_{i\neq j}\left\{ L(W_i)  \beta'(q^4(X_i))L(W_j)  \beta'(q^4(X_j)) T_2 \right\}=o(h^8_{2,n})+o\left ( \frac{1}{n h_{2,n}} \right ). 
\end{eqnarray*}
If $T_1$ is replaced by $T_3$ in (\ref{105tt}), similar to the estimate of $J_3$ in Lemma \ref{Lemma1}, the result is $o(h^{10}_{2,n})$. Together with (\ref{102tt})-(\ref{106tt}), we have
\begin{eqnarray*}
\frac{1}{n^2 h^2_{2,n}}E\sum_{i\neq j} \hat{\lambda}_i Y_i \hat{\lambda}_j =o(h^8_{2,n})+o\left ( \frac{1}{n h_{2,n}} \right ).
\end{eqnarray*}
Similarly, the other terms in (\ref{101tt}) are also $o(h^8_{2,n})+o(1/(nh_{2,n}))$. Hence
\begin{eqnarray*} 
 \frac{1}{n^2 h^2_{2,n}}E\sum_{i\neq j}\hat{\theta}_iY_i \hat{\theta}_j =   o(h^8_{2,n})+o\left ( \frac{1}{n h_{2,n}} \right ).
\end{eqnarray*}
Then by (\ref {90vt}),
\begin{eqnarray}\label {109tt}
 E \left \{\frac{1}{n^2 h^2_{2,n}}\sum_{i=1}^n  \hat{\theta}_iY_i  \sum_{i=1}^n  \hat{\theta}_i \right\}=  o(h^8_{2,n})+o\left ( \frac{1}{n h_{2,n}} \right ) .
\end{eqnarray}
By (\ref{I1a})-(\ref{124Ld}) and (\ref{109tt}),
\begin{eqnarray*}
M_1 =   E\left\{ \bar{g}(t; h_{2,n}) \bar{f}(t; h_{2,n})\right\}- E \bar{g}(t; h_{2,n}) E \bar{f}(t; h_{2,n})+o(h^4_{2,n})+ o\left ( \frac{1}{n h_{2,n}} \right ). 
\end{eqnarray*}
Then by (\ref{I_1}) and (\ref{30cb}), 
\begin{eqnarray}\label {128vbaw}
M_1= \frac{ \sqrt{q(t)}g(t)\mu_{0,2}}{n h_{2,n}} +o(h^4_{2,n}) +o\left ( \frac{1}{n h_{2,n}} \right).
\end{eqnarray}
\textbf{Step 3}. We estimate $M_2$ in (\ref{I2}). Let $\gamma(y)=(1+y)^{-1}$. By Taylor expansion,
\begin{eqnarray*} 
\frac{1}{\hat{f}(t; h_{1,n}, h_{2,n})}= \frac{1}{E\hat{f}(t; h_{1,n}, h_{2,n}) }\left ( 1 + \gamma'(\hat{\rho}_t)\left (\frac{\hat{f}(t; h_{1,n}, h_{2,n})}{E \hat{f}(t; h_{1,n}, h_{2,n})}  -1\right ) \right)
\end{eqnarray*}
where $\hat{\rho}_t$ is between $0$ and $\hat{f}(t; h_{1,n}, h_{2,n})/E\hat{f}(t; h_{1,n}, h_{2,n}) -1$. Hence
\begin{eqnarray}\label {122cb}
&&M_2 = \frac{1}{E \hat{f}(t; h_{1,n}, h_{2,n})}E \left \{\hat{g}(t; h_{1,n}, h_{2,n})\left(\hat{f}(t; h_{1,n}, h_{2,n})- E \hat{f}(t; h_{1,n}, h_{2,n})\right)^2\right\}\nonumber \\
&&+\frac{1}{(E \hat{f}(t; h_{1,n}, h_{2,n}))^2}E \left \{\gamma'(\hat{\rho}_t) \hat{g}(t; h_{1,n}, h_{2,n})  \left (\hat{f}(t; h_{1,n}, h_{2,n})- E \hat{f}(t; h_{1,n}, h_{2,n})\right)^3\right\} \nonumber\\
&&:= \frac{1}{E \hat{f}(t; h_{1,n}, h_{2,n})}M_{2,1}+\frac{1}{(E \hat{f}(t; h_{1,n}, h_{2,n}))^2}M_{2,2}.
\end{eqnarray}
By (\ref {94Eg}) and (\ref{123vbaw}), 
\begin{eqnarray*}\label {101de}
M_{2,1}=M_{2,1,1}+M_{2,1,2}+M_{2,1,3} +M_{2,1,4}
\end{eqnarray*}
where
\begin{eqnarray}\label {99cu}
M_{2,1,1} &=&  E\left\{\bar{g}(t; h_{2,n})  \left ( \bar{f}(t; h_{2,n}) -E\ \bar{f}(t; h_{2,n}) \right )^2\right\}   \nonumber\\ 
M_{2,1,2}&=& 2 E \left \{\bar{g}(t; h_{2,n}) \left( \bar{f}(t; h_{2,n}) -E  \bar{f}(t; h_{2,n}) \right)\frac{1}{n h_{2,n}}\sum_{i=1}^n (\hat{\theta}_i-E \hat{\theta}_i) \right\}\nonumber\\
M_{2,1,3}&=& E\left\{ \bar{g}(t; h_{2,n}) \left ( \frac{1}{n h_{2,n}}\sum_{i=1}^n (\hat{\theta}_i-E \hat{\theta}_i)\right)^2\right\}\nonumber\\
M_{2,1,4}&=& E\left \{ \frac{1}{n h_{2,n}}\sum_{i=1}^n \hat{\theta}_i Y_i \left( \left( \bar{f}(t; h_{2,n}) -E  \bar{f}(t; h_{2,n}) \right)+\frac{1}{n h_{2,n}}\sum_{i=1}^n (\hat{\theta}_i-E \hat{\theta}_i)\right )^2\right\} . \nonumber
\end{eqnarray}
Note that $M_{2,1,1}$ is $I_{2,1}$ in (\ref{40vbb}) with $h_n$ replaced by $h_{2,n}$. Then  by (\ref{37cba}),
\begin{eqnarray}\label {M211}
M_{2,1,1} = \frac{\sqrt{ q(t)}g(t)f(t)\mu_{0,2}}{n h_{2,n}} +o\left(\frac{1}{n h_{2,n}}\right).
\end{eqnarray}
By Lemma \ref{Lemma1} and Lemma \ref{Lemma2}, and similar to (\ref{92pa}), 
$$E (  E  \{ (\hat{\theta}_i-E \hat{\theta}_i)|X_i,  X_j, X_k   \}| X_j, X_k)=o(h^4_{2,n}).$$ 
Denote $F_i = K(W_i)\alpha(q(X_i)).$ 
Then
\begin{eqnarray}\label {106pa}
M_{2,1,2} &=& \frac{2}{n^3 h^3_{2,n}} \sum_{i , j, k \ are \ not \ equal} E\bigg \{F_k Y_k  \left \{F_j -E F_j \right \}\times E \left \{ E   (\hat{\theta}_i-E \hat{\theta}_i)| X_j, X_k \right \}\bigg\} \nonumber\\&& + o\left (  \frac{1}{n h_{2,n}}\right )\nonumber\\
&=& \frac{o(h^4_{2,n})}{h_{2,n}}  E \left |F_k Y_k \right| \frac{1}{h_{2,n}} E\left |F_j-EF_j \right | + o\left (  \frac{1}{n h_{2,n}}\right )\nonumber\\
&=&o(h^4_{2,n})+ o\left (  \frac{1}{n h_{2,n}}\right ).
\end{eqnarray}
Similarly,
\begin{eqnarray}\label {M213}
M_{2,1,3}+M_{2,1,4}=o(h^4_{2,n})+ o\left (  \frac{1}{n h_{2,n}}\right ).
\end{eqnarray}
By (\ref{M211})-(\ref{M213}),
\begin{eqnarray}\label {109deh}
 M_{2,1}=   \frac{\sqrt{ q(t)}g(t)f(t)\mu_{0,2}}{n h_{2,n}}+  o(h^4_{2,n})+ o\left (  \frac{1}{n h_{2,n}}\right ).
\end{eqnarray}
For $M_{2,2}$ defined in (\ref{122cb}), by H\a"{o}lder's inequality,\\
\begin{eqnarray}\label {124tt}
|M_{2,2}| &\leq& \|\gamma'(\cdot)\|_{\infty}  \left \{E  \left(   \hat{g}(t; h_{1,n}, h_{2,n})\right)^2\right\}^{1/2} \nonumber\\
&& \times \left \{E \left ( \bar{f}(t; h_{2,n}-E \bar{f}(t; h_{2,n})+\frac{1}{n h_{2,n}}\sum_{i=1}^n (\hat{\theta}_i-E \hat{\theta}_i)\right)^6\right \}^{1/2}. \qquad
\end{eqnarray}
By (\ref{Ebarf}), 
\begin{eqnarray}\label {99vt}
 E \left ( \bar{f}(t; h_{2,n}-E \bar{f}(t; h_{2,n})   \right )^6= O\left (\frac{1}{n^3 h^3_{2,n}}  \right).
\end{eqnarray}
Similar to the estimate of (\ref{90vt}), \\
\begin{eqnarray}\label {100vt}
E \left ( \frac{1}{n h_{2,n}}\sum_{i=1}^n (\hat{\theta}_i-E \hat{\theta}_i)\right)^6=  o(h^{8}_{2,n})+ o\left ( \frac{1}{n^2 h^2_{2,n}} \right).
\end{eqnarray}
Applying (\ref{99vt}) and (\ref{100vt}) to (\ref{124tt}),  
\begin{eqnarray}\label{112de}
M_{2,2}= o (h^4_{2,n}) + o\left ( \frac{1}{n h_{2,n}} \right ). 
\end{eqnarray}
Since $E \bar{f}(t; h_{1,n}, h_{2,n})=f(t)+o(1)$, by (\ref{122cb}), (\ref{109deh}) and (\ref{112de}),  
\begin{eqnarray}\label {126cb}
M_2=\frac{\sqrt{q(t)}g(t)\mu_{0,2}}{n h_n} +o(h^4_{2,n})+ o\left ( \frac{1}{n h_{2,n}} \right ).
\end{eqnarray}
By (\ref{75vbt}), (\ref{96vbt}), (\ref{128vbaw}) and (\ref{126cb}), we obtain (\ref{Thm2(1)}). \\
Now we prove (\ref{Thm2(2)}). By the Taylor expansion of $1/\hat{f}(t; h_{1,n}, h_{2,n})$ at the beginning of Step 3,
\begin{eqnarray}\label {153svb}
&&   \hat{r}(t; h_{1,n}, h_{2,n})- r(t)\nonumber\\
&=&    \frac{f(t)\hat{g}(t; h_{1,n}, h_{2,n})- g(t)\hat{f}(t; h_{1,n}, h_{2,n}) }{f(t)\hat{f}(t; h_{1,n}, h_{2,n})}  \nonumber\\
&=& \frac{N_1}{f(t) E\hat{f}(t; h_{1,n}, h_{2,n})  }+ \frac{N_2}{f(t) (E\hat{f}(t; h_{1,n}, h_{2,n}))^2  } 
\end{eqnarray}
where
\begin{eqnarray}\label {J1}
N_1 &=&   f(t)\hat{g}(t; h_{1,n}, h_{2,n})- g(t)\hat{f}(t; h_{1,n}, h_{2,n})   \nonumber\\
N_2 &=&   \left \{  f(t)\hat{g}(t; h_{1,n}, h_{2,n})- g(t)\hat{f}(t; h_{1,n}, h_{2,n})  \right \}\nonumber\\
&& \times \left \{ \gamma'(\hat{\rho}_t) \left (  \hat{f}(t; h_{1,n}, h_{2,n}) -E\hat{f}(t; h_{1,n}, h_{2,n}) \right )\right \}. \label {J2}
\end{eqnarray}
Since  $E\hat{f}(t;  h_{1,n}, h_{2,n})=f(t)(1+o(1))$, then
\begin{eqnarray}\label {127et}
E \left \{\hat{r}(t; h_{1,n}, h_{2,n})- r(t)  \right \}^2= \left (\frac{E N_1^2}{f^4(t) }+  \frac{2 E (N_1 N_2)}{f^5(t) } + \frac{E N_2^2}{f^6(t) }\right)(1+o(1)).
\end{eqnarray}
By (\ref{94Eg}) and (\ref{123vbaw}), 
\begin{eqnarray}\label {117ua}
E N^2_1&=&  E \left \{  f(t) \bar{g}(t; h_{2,n}) - g(t)\bar{f}(t; h_{2,n})    \right \}^2 + E \left \{ \frac{1}{n h_{2,n}}\sum_{i=1}^n  \hat{\theta}_i \left \{f(t)Y_i -g(t)\right\} \right \}^2\nonumber\\
&& +2  E\left \{ f(t) \bar{g}(t; h_{2,n})   - g(t) \bar{f}(t; h_{2,n}) \right \} \left \{ \frac{1}{n h_{2,n}}\sum_{i=1}^n  \hat{\theta}_i\left \{f(t)Y_i -g(t)\right\} \right \} \nonumber\\
&:=&E N_{1,1}+E N_{1,2}+2 E N_{1,3}.
\end{eqnarray}
By (\ref{55svb}),
\begin{eqnarray} \label {120eg}
E N_{1,1}=  \frac{f^{7/2}(t)|r'(t)|^{1/4}  \sigma^2(t) \mu_{0,2}}{n h_n} +f^4(t)\theta^2(t)h^8_n +o(h^8_n)+o\left ( \frac{1}{nh_n}\right).
\end{eqnarray}
Similar to the estimate of (\ref{90vt}) above, 
\begin{eqnarray}\label {121eg}
E N_{1,2}= o(h^8_{2,n})+o\left (\frac{1}{n h_{2,n}}  \right ).
\end{eqnarray}
By H\a"{o}lder's inequality, (\ref{120eg}) and (\ref{121eg}),
\begin{eqnarray}\label {120ua}
|E N_{1,3}| \leq (E N_{1,1})^{1/2}(E N_{1,2})^{1/2} = o(h^8_{2,n})+ o\left (\frac{1}{n h_{2,n}}  \right ).
\end{eqnarray}
By (\ref{117ua})-(\ref{120ua}),
\begin{eqnarray}\label {J1a}
E N^2_1=\frac{f^{7/2}(t)|r'(t)|^{1/4}  \sigma^2(t) \mu_{0,2}}{n h_n} +f^4(t)\theta^2(t)h^8_n +o(h^8_n)+o\left ( \frac{1}{nh_n}\right).
\end{eqnarray}
Next we estimate $EN^2_2$. By H\a"{o}lder's inequality, 
\begin{eqnarray}\label {130vw}
E N^2_{2}\leq \|\gamma'(\cdot)\|^2_{\infty} \left ( E N_1^4 \right )^{1/2} \left ( E \left (  \hat{f}(t; h_{1,n}, h_{2,n}) -E\hat{f}(t; h_{1,n}, h_{2,n}) \right )^4  \right )^{1/2}.
\end{eqnarray}
By (\ref{94Eg}), (\ref{123vbaw}), (\ref{53wb}) and similar to the estimate of (\ref{90vt}), 
\begin{eqnarray}\label {136vk}
E N_1^4 &=&E \left \{f(t) \bar{g}(t; h_{2,n})  -g(t) \bar{f}(t; h_{2,n} )+  \frac{1}{n h_{2,n}}\sum_{i=1}^n  \hat{\theta}_i \left \{f(t)Y_i -g(t)\right\} \right \}^4\nonumber\\
&\leq&8 E \left \{ f(t) \bar{g}(t; h_{2,n} ) - g(t) \bar{f}(t; h_{2,n} )  \right \}^4 +8 E \left \{    \frac{1}{n h_{2,n}}\sum_{i=1}^n  \hat{\theta}_i \left \{f(t)Y_i -g(t)\right\} \right \}^4\nonumber\\
&=& o(h^8_{2,n})+ o\left ( \frac{1}{n h_{2,n}} \right ).
\end{eqnarray}
By (\ref{123vbaw}), (\ref{54wb}) and similar to the estimate of (\ref{90vt}), 
\begin{eqnarray}\label {137vk}
&&E \left (  \hat{f}(t; h_{1,n}, h_{2,n}) -E\hat{f}(t; h_{1,n}, h_{2,n}) \right )^4\nonumber\\
 &\leq& 8 E \left \{\bar{f}(t; h_{2,n} )-E \bar{f}(t; h_{2,n} )\right \}^4+ 8 E \left \{ \frac{1}{n h_{2,n}}\sum_{i=1}^n (\hat{\theta}_i-E \hat{\theta}_i)\right \}^4 \nonumber\\
&=&  o(h^8_{2,n})+ o\left ( \frac{1}{n h_{2,n}}\right).
\end{eqnarray}
Applying (\ref{136vk}) and (\ref{137vk}) to (\ref{130vw}), we have
\begin{eqnarray}\label {127eh}
E N^2_2=  o(h^8_{2,n})+ o\left (\frac{1}{n h_{2,n}}  \right ).
\end{eqnarray}
Since $|E (N_1 N_2)| \leq (E N_1^2)^{1/2}(E N_2^2)^{1/2}$, then applying (\ref{J1a}) and (\ref{127eh}) to (\ref{127et}), we obtain (\ref{Thm2(2)}). The integrated mean squared error\\
\begin{eqnarray*}
&&\int_{t \in D_{rf}}  E \left ( \hat{r}(t; h_{1,n}, h_{2,n})- r(t) \right )^2 dt\nonumber\\
&&=h^8_{2,n}  \int_{t \in D_{rf}}  \theta^2(t)dt + \frac{\mu_{0,2}\sigma^2}{n h_{2,n}} \int_{t \in D_{rf}}  \frac{|r'(t)|^{1/4}}{\sqrt{f(t)}}dt + o(h^8_{2,n})+o\left (\frac{1}{n h_{2,n}}  \right ).\qquad
\end{eqnarray*}
Taking derivative with respect to $h_{2,n}$ and letting the result equal to $0$, we obtain the optimal bandwidth 
\begin{eqnarray*}
h_n^*=\left (\frac{1}{n}  \right )^{1/9}\left ( \frac{\mu_{0,2}\sigma^2  }{8 \int_{t \in D_{rf}} \theta^2(t)dt} \int_{t \in D_{rf}} \frac{|r'(t)|^{1/4}}{\sqrt{f(t)}}dt\right )^{1/9}.
\end{eqnarray*}
\end{proof}

\begin{proof}[Proof of Theorem \ref{central2}]
The proof is similar to the proof of Theorem \ref{central1} and use the results in Theorem \ref{Thm2}.
\end{proof}




\begin{thebibliography}{99}

\bibitem{Abramson}
I. Abramson.
On bandwidth variation in kernel estimates - a square-root law.
{\it Ann. Statist.} {\bf 10} (1982), 1217-1223.

\bibitem{HJB1983}
H. J. Bierens. Uniform consistency of kernel estimators of a regression function
under generalized conditions. {\it Journal of the American Statistical Association} {\bf 79} (1983), 699-707.

\bibitem{GC1977} 
G. Collomb. Estimation non param\'{e}trique de la r\'{e}gression par la m\'{e}thode du noyau: propri\'{e}t\'{e}  de convergence asymptotiquement normale ind\'{e}pendante Annales de le Facult\'{e}des  sciences de I'Universit\'{e} de Clermont. S\'{e}rie Math\'{e}matiques des sciences de
I'Universit\'{e} de Clermont. S\'{e}rie  Math\'{e}matiques {\bf 65} (1977), 24-46.

\bibitem{EM2000}
U. Einmahl and D. M. Mason. An empirical process approach to the uniform consistency
of kernel-type function estimators, {\it Journal of Theoretical Probability} {\bf 13} (2000), 1-37.

\bibitem{EM2005}
U. Einmahl and D. M. Mason.
Uniform in bandwidth consistency of kernel-type function estimators, {\it Ann.
Statist.} {\bf 33} (2005), 1380-1403.


\bibitem{GM1984}
T. Gasser and H. G. M\a"{u}ller. Estimating regression functions and their derivatives by the kernel method. {\it Scandanavian Journal of Statistics.} \textbf{11} (1984), 171-185.


\bibitem{GineSang}
E. Gin\'{e} and H. Sang.
Uniform asymptotics for kernel density estimators with variable bandwidths.
{\it J. Nonparametr. Stat.} {\bf 22} (2010), 773-795.

\bibitem{GS2013}
E. Gin\'{e} and H. Sang.
On the estimation of smooth densities by strict probability densities at optimal rates in sup-norm. 
{\it IMS Collections, From Probability to Statistics and Back: High-Dimensional Models and Processes}
{\bf 9} (2013), 128-149.

\bibitem{Hall}
P. Hall. On the bias of  variable bandwidth kernel estimators. {\it Biometrika} {\bf 77} (1990), 529-535.

\bibitem{HallHuMarron}
P. Hall, T. Hu and J. S. Marron.
Improved Variable Window Kernel Estimates of Probability Densities.  {\it Ann. Statist.} {\bf 23} (1995), 1-10. 

\bibitem{Hay}
Tristen Hayfield and Jeffrey S. Racine.
Nonparametric Econometrics: The np Package. Journal of Statistical Software.  {\it Ann. Statist.} {\bf 25} (2008), 5.
URL http://www.jstatsoft.org/v27/i05/.

\bibitem{HM1988}
P. Hall and J. S. Marron.
Variable Window Width Kernel Estimates of Probability Densities. {\it Probab. Theory Related Fields} {\bf 80}  (1988), 37-49. Erratum: {\it Probab. Theory Related Fields} {\bf 91} 133.


\bibitem{JonesMcKayHu}
M. C. Jones, I. J. McKay and T.-C. Hu. Variable location and scale kernel density estimation. {\it Ann. Inst. Statist. Math.}  {\bf 46} (1994), 521-535.


\bibitem{McKay a}
I. J. McKay. A note on bias reduction in variable kernel density estimates. {\it Canad. J. Statist.} {\bf 21} (1993a), 367-375.

\bibitem{McKay b}
I. J. McKay. Variable kernel methods in density estimation.  Ph.D Dissertation, Queen's University, (1993b).


\bibitem{MS}
H. M\"uller and U. Stadtm\"uller. Variable bandwidth kernel estimators of regression curves,
{\it Ann. Statist.} {\bf 15} (1987), 182-201.


\bibitem{EAN1964}
E. A. Nadaraya. On estimating regression. {\it Theory of Probability and Its Applications}
{\bf 9} (1964), 141-142.


\bibitem{NS2016}
J. Nakarmi and H. Sang. Central limit theorem for the variable bandwidth kernel density estimators,  {\it Journal of the Korean Statistical Society}, 47 (2018), no. 2, 201-215.

%

\bibitem{KN1976}
K. Noda. Estimation of a regression function by the Parzen kernel-type density estimators, {\it Annals of the Institute of Statistical Mathematics} {\bf 28} (1976), 221-234.

\bibitem{R}
M. Rosenblatt. Conditional Probability Density and Regression Estimators.
{\it Multivariate Analysis II} (1969),  Academic Press, New York, 25-31.

\bibitem{Novak}
S. Yu. Novak. A generalized kernel density estimator. (Russian) {\it Teor. Veroyatnost. i Primenen.} {\bf 44} (1999) 634--645; translation in {Theory Probab. Appl.} {\bf 44} (2000) 570-583.


\bibitem{MS2000}
M. G. Schimek.  {\it Smoothing and regression: Approaches, computation, and application}. John Wiley \& Sons, (2000).

%
%


\bibitem{TerrellScott}
G. R. Terrell and D. Scott. Variable kernel density estimation. {\it Ann. Statist.} {\bf 20} (1992), 1236-1265.

\bibitem{wasserman}
L. Wasserman. {\it All of nonparametric statistics}, Springer, (2006).

\bibitem{GSW1964}
G. S. Watson. Smooth regression analysis, {\it Sankhya}. Series A {\bf 26} (1964), 359-372.


\end{thebibliography}
\end{document}